\documentclass[12pt]{article}
\usepackage{graphicx}
\usepackage[bf]{caption2}
\usepackage{amsfonts, amssymb, amsthm, amsmath}
\usepackage{enumerate}
\newcounter{theorem}

\newcounter{theoremcounter}

\newcounter{corollarycounter}
\newtheorem{theorem}[theoremcounter]{Theorem}

\newtheorem{corollary}[corollarycounter]{Corollary}
\newcommand{\la}{\lambda}
\newcommand{\vp}{{\mathbf p}}
\newcommand{\vz}{{\mathbf z}}
\newcommand{\vf}{{\mathbf f}}

\begin{document}

%%%%% TITLE %%%%%

\title{Bounds on the rate of convergence for inhomogeneous $M/M/S$ systems with either
state-dependent transitions, or batch arrivals and service, or both} %in capital

\date{}
%%%%% AUTHORS %%%%%

\author{Alexander Zeifman \and Anna Korotysheva  \and Yacov Satin    \and Rostislav Razumchik \and Victor Korolev  \and Ksenia Kiseleva}

\maketitle

\textbf{Abstract.}
In this paper one presents method
for the computation of convergence bounds
for four classes of multiserver queueing systems,
described by inhomogeneous Markov chains.
Specifically one considers
inhomogeneous $M/M/S$ queueing system with possibly
state-dependent arrival and service intensities
and additionally possible batch arrivals
and batch service.
The unified approach based on
logarithmic norm of linear operators
for obtaining
sharp upper and lower bounds on the rate of convergence
and corresponding sharp perturbation bounds
is described.
As a side result, one shows by virtue of numerical
examples that
the approach based on logarithmic norm can also be used
for approximation of limiting characteristics (idle probability
and mean number of customers in the system) of the
considered systems with given approximation error.
Extensive numerical examples are provided.

\section{Introduction}

In this paper one considers
the class of Markov processes,
which is usually used to describe the evolution of
the total number of customers
in inhomogeneous Markov queueing
systems.
Suppose that the system's
state space is $\mathcal{X}=\{ 0, 1, 2 \dots \}$,
where the $(i)$-th state means that there are $i$ customers
in the system.
Thought the paper it is assumed that
all possible transition intensities between states are non-random functions of time
and may depend on the state of the system state.

There are two most common problems related to such systems: the computation of
time-dependent distribution of the state probabilities
and the limiting distribution (for example, in case of periodic intensities);
computation of the rate of convergence and perturbation
bounds.

This paper deals with the second problem related to
the following classes of inhomogeneous Markov queueing
systems with arbitrary finite number of servers $S$:
\begin{enumerate}[I.]
  \item inhomogeneous $M/M/S$ queueing system with possibly state-dependent arrival and service intensities;
  \item inhomogeneous $M/M/S$ queueing system with state-independent batch arrivals and state-dependent service intensity;
  \item inhomogeneous $M/M/S$ queueing system with state-dependent arrival intensity and batch service;
  \item inhomogeneous $M/M/S$ queueing system with state-independent batch arrivals and batch service.
\end{enumerate}
\noindent Here for these four system classes
one describes the unified approach based on
logarithmic norm of linear operators
for obtaining
sharp upper and lower bounds on the rate of convergence
and corresponding sharp perturbation bounds.

This unified approach has already been successfully applied to
system from the  ${\rm I}^{st}$ and ${\rm IV}^{th}$ class.
Specifically
for the inhomogeneous $M/M/1$ system with state-dependent arrival
and service intensities, as well as for the state-independent
inhomogeneous $M/M/S$ system
the bounds were firstly obtained in \cite{Zeifman1991},
\cite{Zeifman1995b} and \cite{Doorn2010}.
Systems belonging to the ${\rm IV}^{th}$ class have been studied
in a number of papers (see, for example, \cite{Nelson1988}, \cite{Li2016}, \cite{Li2016-2})
and the results related to convergence have
been also obtained in \cite{Satin2013,Zeifman2014a}.
Here one demonstrates that the approach is also suitable for
the systems from the ${\rm II}^{nd}$ and ${\rm III}^{d}$ class
and thus offers a unified way toward analysis of ergodicity properties
of such Markov chains.

The approach is based on the special
properties of linear systems of differential equations with
non-diagonally non-negative matrices.
Specifically, if the column-wise sums of the
elements of this matrix are identical and equal to, say,
$-\alpha^*(t)$, then the exact upper bound of order
$\exp\big\{-\int_0^t\alpha^*(u)\,du\big\}$ can be obtained for the
rate of convergence of the solutions of the system
in the corresponding metric. Moreover, if the column-wise sums of
the absolute values of the elements of this matrix are identical and
equal to, say, $\chi^*(t)$, then the exact lower bound of order
$\exp\big\{-\int_0^t\chi^*(u)\,du\big\}$ can be obtained for the
convergence rate as well. The bounds are obtained in three steps.
At first step one excludes the $(0)$ state from
the forward Kolmogorov system of differential equations
and thus obtains the new system
with the new intensity matrix which is, in general, not non-diagonally
non-negative. The second step is to transform the new intensity matrix
in such a way that non-diagonally elements are non-negative
and which leads to (loosely speaking) least distance between
specifically defined upper and lower bounds. At third step one uses
the logarithmic norm for the estimation of the convergence rate.

Here the key step is the second one. The transformation
is made using a sequence of positive numbers $\{d_i, i \ge 1\}$,
which does not have any probabilistic sense and
can be considered as an analogue of Lyapunov functions.
For the detailed discussion on application of logarithmic norm and related
techniques one can refer to the series of papers
\cite{Doorn2010,Granovsky2004,Zeifman1995b,Zeifman2006,Zeifman2014a,Zeifman2015s}.

The advantages of this three-step approach
is that it allows one to deal with time-homogeneous
and time-inhomogeneous processes and
it leads to exact both upper and lower bounds for the
convergence rate.
In time-homogeneous case (of the four classes of systems introduced above)
the approach allows one to obtain the
correspondent bounds for the decay parameter and gives an
explicit bounds in total variation norm (see \textit{Theorem 2}).

The proposed approach  allows one also to address the problem
of computation of the limiting distribution of the inhomogeneous
 Markov chain from a different perspective.
In general there are several approaches, which allow one to
obtain more or less accurate solutions.
These are the exact and approximate
numerical solution of the system of differential equations,
approaches assuming piecewise constant parameters
and approaches based on modified system characteristics.
For the review of many results one can refer to \cite{Li2016-4}.
Although using the proposed approach one cannot
determine the state probabilities as functions of time $t$,
it is possible to compute approximately the limiting distribution,
while having analytically computable expressions for the approximation errors.
Using truncation techniques, which were developed in
\cite{Zeifman2006,Zeifman2014i}, in the numerical
one presents the results of the computation of the
limiting characteristics in inhomogeneous $M/M/S$ systems
of each of the four classes described above.
The most interesting insight from the experiments is the following.
Choose the arrival and service intensities
in inhomogeneous $M/M/S$ (${\rm I}^{st}$ class).
Then if one uses these intensities (after a certain
modification allowing
bulk arrivals and group services) in inhomogeneous $M/M/S$
from the ${\rm II}^{nd}$, ${\rm III}^{d}$
or ${\rm IV}^{th}$ class, then the
limiting mean number of customers for both systems
coincide,
while idle probabilities do not.

The paper is structured as follows.
In the next section one gives the general description
of the system under consideration
and introduces the necessary notation.
Section 3 contains the main result of the paper
i.e. the theorem which specifies the
convergence bounds.
Section 4 provides explicit expressions
for functions needed to compute
convergence bounds for four special cases
of the considered system.
In the last two sections one provides
extensive
numerical examples and gives directions of
further research.

\section{System description and definitions}

Let the integer-valued time-dependent random variable  $X(t)$ denote the
total number of customers at time $t$ in a markovian queueing system.
Then the process $\{ X(t), \ t\geq 0\}$ is
a (possibly inhomogeneous) continuous-time Markov
 chain with state space $\mathcal{X}=\{ 0, 1, 2 \dots \}$.
Denote by $p_{ij}(s,t)=P\left\{ X(t)=j\left| X(s)=i\right. \right\}$,
$i,j \ge 0, \;0\leq s\leq t$ the transition probabilities of
$X(t)$ and by  $p_i(t)=P \left\{ X(t) =i \right\}$ -- the
probability  that the Markov chain $X(t)$ is in state $i$ at time $t$.
Let $\vp(t) = \left(p_0(t), p_1(t), \dots\right)^T$ be
probability distribution vector at instant $t$.
Throughout the paper we assume that in an element of time $h$ the
possible transitions and their associated probabilities are
\begin{equation}
p_{ij}(t,t+h)=
%\Pr\left(X\left( t+h\right) =j/X\left( t\right) =i\right) =
\left\{
\begin{array}{cc}
q_{ij}(t)  h+\alpha_{ij}\left(t, h\right), & \mbox {if
}j\neq i
\\ 1+ q_{ii}(t) h+\alpha_{i}\left(
t,h\right), & \mbox {if } j=i,
\end{array}
\right.
i, j \in \mathcal{X},
  \label{4001}
\end{equation}
\noindent where all  $\alpha_{i}(t,h)$ are $o(h)$ uniformly in $i$,
i.e.  $\sup_i |\alpha_i(t,h)| = o(h)$
and
$$
q_{ii}(t) = - \sum_{k \in \mathcal{X}, k\neq i}q_{ik}(t).
$$

Applying the standard approach developed in
\cite{Granovsky2004,Zeifman1995b,Zeifman2006} it is assumed that
that all the intensity functions $q_{ij}(t)$ are linear combinations of a
finite number of locally integrable on $[0,\infty)$ non-negative
functions.

The matrix $Q(t)=(q_{ij}(t))_{i,j=0}^{\infty}$ is the
intensity matrix of the chain $\{ X(t), \ t\geq 0\}$.
Henceforth it is assumed that the $Q(t)$ is
essentially bounded, i. e.
\begin{equation}
\sup_i|q_{ii}(t)| = L(t) \le L < \infty, \label{0102-1}
\end{equation}
\noindent for almost all $t \ge 0$.

Probabilistic dynamics of
the process $\{ X(t), \ t\geq 0\}$
is given by the
forward Kolmogorov system
\begin{equation} \label{ur01}
\frac{d}{dt}\vp(t)=A(t)\vp(t),
\end{equation}
\noindent where $A(t) = Q^T(t)$ is the {transposed} intensity
matrix.
%and $\vp(t)$ is the column vector of state probabilities,
%$\vp(t) = \left(p_0(t), p_1(t), \dots\right)^T$.

Throughout the paper by $\|\cdot\|$  we denote  the $l_1$-norm, i.
e.  $\|{\vp(t)}\|=\sum_{k\in \mathcal{X}} |p_k(t)|$, and
$\|Q(t)\| = \sup_{j \in \mathcal{X}} \sum_{i\in \mathcal{X}} |q_{ij}|$.
Let $\Omega$ be a set all stochastic vectors, i. e. $l_1$ vectors
with non-negative coordinates and unit norm. Hence we have $\|A(t)\|
= 2\sup_{k\in \mathcal{X}} |q_{kk}(t)| \le 2 L $ for almost all $t \ge 0$.
Hence the operator function $A(t)$ from $l_1$ into itself is bounded
for almost all $t \ge 0$ and locally integrable on $[0;\infty)$.
Therefore we can consider (\ref{ur01}) as a differential equation in
the space $l_1$ with bounded operator.

It is well known (see \cite{Daleckij1974}) that the Cauchy problem
for differential equation (\ref{ur01}) has a unique solutions for an
arbitrary initial condition, and  $\vp(s) \in \Omega$ implies
$\vp(t) \in \Omega$ for $t \ge s \ge 0$.

Denote by $E(t,k) = E(X(t)|X(0)=k)$ the conditional expected
number of customers in the system at instant $t$, provided that
initially (at instant $t=0$) $k$ customers were present in the system.
Then $E_{p}(t)=\sum_{k\ge 0}E(t,k)p_k(0)$ is the unconditional
expected number of customers in the system at instant $t$, given that
the initial distribution of the total number of customers was
$\vp(0)$.

In order to obtain perturbation bounds we consider a class of
perturbed Markov chains  $\{\bar{X}(t), t\geq 0\}$
defined on the same state space $\mathcal{X}$ as the original
Markov chain $\{{X}(t), t\geq 0\}$,
with the intensity matrix $\bar{A}(t)$ and the same
restrictions as imposed on $A(t)$. It is assumed that $\| \hat{A}(t)\|=
\|A(t)-\bar{A}(t)\|
\le \varepsilon$, for almost all $t \ge 0$,
which means the perturbations are considered to be small.

Before proceeding to the derivation of the main results of the paper,
we recall two definitions.
Recall that a Markov chain $\{ X(t), \ t\geq 0\}$ is called {\it weakly ergodic}, if
$\|\vp^{*}(t)-\vp^{**}(t)\| \to 0$ as $t \to \infty$ for any
initial conditions $\vp^{*}(0)$ and $\vp^{**}(0)$, where $\vp^{*}(t)$ and $\vp^{**}(t)$ are the corresponding solutions of
(\ref{ur01}).
A Markov chain $\{ X(t), \ t\geq 0\}$ has  the {limiting mean} $\varphi (t)$, if
$ \lim_{t \to \infty }  \left(\varphi (t) - E(t,k)\right) = 0$ for
any $k$.

\section{Main results}

Recall that one has introduced $A(t)$ as the transposed intensity matrix $Q(t)$.
Thus it has the form
\begin{equation}
A(t)=\left(
\begin{array}{cccccc}
a_{00}(t) & a_{01}(t)  &   \cdots & a_{0r}(t) & \cdots\\
a_{10}(t)  & a_{11}(t)  &  \cdots & a_{1r}(t) & \cdots\\
a_{20}(t)  & a_{21}(t)    & \cdots & a_{2r}(t) & \cdots\\
\cdots \\
a_{r0}(t) & a_{r1}(t) &  \cdots    &  a_{rr}(t) & \cdots\\
\cdots
\end{array}
\right),
\end{equation}
%}
\noindent where  $a_{ii}(t)=-\sum_{k \in \mathcal{X}, k \neq i} a_{ki}(t)$.
Since $p_0(t) = 1 - \sum_{i = 1}^\infty p_i(t)$ due to
normalization condition, one can rewrite the system (\ref{ur01})
as follows:
\begin{equation}
\frac{d}{dt}\vz(t)= B(t)\vz(t)+\vf(t), \label{2.06}
\end{equation}
\noindent where
$$
\vf(t)=\left( a_{10}(t),  a_{20}(t),\dots \right)^{T}, \
\vz(t)=\left(p_1(t), p_2(t),\dots \right)^{T},
$$
\begin{equation}
{\footnotesize
B \!=\! \left(b_{ij}(t)\right)_{i,j=1}^{\infty} \!=\! \left(
\begin{array}{ccccc}
a_{11}\left( t\right) \!-\!a_{10}\left( t\right) & a_{12}\left( t\right)
\!-\!a_{10}\left( t\right) & \cdots & a_{1r}\left( t\right)
\!-\!a_{10}\left(
t\right) & \cdots \\
a_{21}\left( t\right) \!-\!a_{20}\left( t\right) & a_{22}\left( t\right)
\!-\!a_{20}\left( t\right) & \cdots & a_{2r}\left( t\right)
\!-\!a_{20}\left(
t\right) & \cdots \\
\cdots & \cdots & \cdots & \cdots  & \cdots \\
a_{r1}\left( t\right) \!-\!a_{r0}\left( t\right) & a_{r2}\left( t\right)
\!-\!a_{r0}\left( t\right) & \cdots & a_{rr}\left( t\right)
\!-\!a_{r0}\left( t\right) & \cdots \\
\vdots & \vdots & \vdots & \vdots  & \ddots
\end{array}
\right).}\label{2.07}
\end{equation}
See detailed discussion of this transformation in
\cite{Granovsky2004,Zeifman1995b,Zeifman2006}.
Let $\{d_i, \ i \ge 1\}$
with $d_1 = 1$  be an increasing sequence of positive numbers. Put
\begin{equation}
W=\inf_{i \ge 1} \frac {d_i}{i}. %, \quad g_i=\sum_{n=1}^i d_n, \, d=\inf_{i \ge 1} d_i.
\label{2002}
\end{equation}
and denote by $D$ the upper triangular matrix of the following form:
\begin{equation}
D=\left(
\begin{array}{ccccccc}
d_1   & d_1 & d_1 & \cdots  \\
0   & d_2  & d_2  &   \cdots  \\
0   & 0  & d_3  &   \cdots  \\
\vdots & \vdots & \vdots & \ddots \\
\end{array}
\right). \label{2013}
\end{equation} \noindent
Let $l_{1D}$ be
the corresponding space of sequences
$$l_{1D}=\left\{{\bf z}(t) = (p_1(t),p_2(t),\cdots)^{T} |\, \|{\bf z}(t)\|_{1D} \equiv \|D {\bf z}(t)\|_1 <\infty \right\}
$$
and introduce also the auxiliary norm $\|\cdot\|_{1E}$ defined as
$\|{\bf z}(t)\|_{1E}=\sum_{k=1}^\infty k|p_k(t)|$.
%if ${\bf x} = (x_1,x_2,\cdots)^{T}$.
Then in $\|\cdot\|_{1D}$ norm the following two inequalities
 hold:
\begin{eqnarray}
\|{\bf z}(t)\|_{1D} = d_1 \left|\sum\limits_{i=1}^{\infty} p_i(t) \right|
+ d_2 \left|\sum\limits_{i=2}^{\infty} p_i(t) \right| + d_3
\left|\sum\limits_{i=3}^{\infty} p_i(t) \right| + \dots \ge \nonumber
\\ \ge  \left(\left|\sum\limits_{i=1}^{\infty} p_i(t) \right| +
\left|\sum\limits_{i=2}^{\infty} p_i(t) \right| +
\left|\sum\limits_{i=3}^{\infty} p_i(t) \right| + \dots \right) \ge
\nonumber \\ \ge \frac{1}{2}
\left(\left(\left|\sum\limits_{i=1}^{\infty} p_i(t) \right| +
\left|\sum\limits_{i=2}^{\infty} p_i(t) \right|\right) +
\left(\left|\sum\limits_{i=2}^{\infty} p_i(t) \right| +
\left|\sum\limits_{i=3}^{\infty} p_i(t) \right|\right) + \dots \right)
\ge  \nonumber \\ \ge \frac{1}{2} \sum\limits_{i=1}^{\infty} |p_i(t) |
= \frac{1}{2} \|{\bf z}(t)\|_1, \label{normD}
\end{eqnarray}
\begin{small}
\begin{eqnarray} \label{normDE}
\|\vz(t)\|_{1E}
=
\sum_{k=1}^{\infty} k|p_k(t)|
=
\sum_{k=1}^{\infty}\frac{k}{d_k}d_k|p_k(t)|
\le W^{-1}\sum_{k=1}^{\infty}d_k|p_k(t)|
=
\nonumber
\\
=
W^{-1}\sum_{k=1}^{\infty} d_k \left| \sum_{i= k}^\infty p_i(t) -\sum_{i= k-1}^\infty p_i(t) \right|
\le W^{-1}\sum_{k=1}^{\infty} d_k \left(\left| \sum_{i= k}^\infty p_i(t)\right|
+\left| \sum_{i= k-1}^\infty p_i(t) \right| \right)  \le
\nonumber
\\
\le \frac{2}{W}\sum_{k=1}^{\infty} d_k \left| \sum_{i= k}^\infty p_i(t)\right| \le
\frac{2}{W}\| \vz(t)\|_{1D}.
\end{eqnarray}
\end{small}

Consider the equation (\ref{2.06}) in the space $l_{1D}$, where $B(t)$
and ${\bf f}(t)$ are locally integrable on $[0, +\infty)$.
Let one compute the logarithmic norm of operator function ${B}(t)$.
The motivation behind this can be found in \cite{Doorn2010} and detailed proofs
are provided in \cite{Zeifman1995a}. Recall that the logarithmic norm of operator
function ${B}(t)$ is defined as
$$\gamma({B}(t)) = \lim_{h \to
+0}h^{-1}\left(\|I+hB(t)\|-1\right).$$
\noindent Denote by $V(t, s) = V(t)V^{-1}(s)$ the Cauchy operator of the equation (\ref{2.06}).
Then the  important inequality holds
$$e^{-\int_s^{t} \gamma(-B(u))\, du} \le \|V(t, s)\| \le e^{\int_s^{t} \gamma(B(u))\, du}.$$
Further, for an operator function from $l_1$ to itself
one has the simple formula
$$
\gamma({B}(t)) = \sup_j
\left(b_{jj}(t)+\sum_{i \neq j} |b_{ij}(t)|\right).
$$
Moreover for the logarithmic norm
of the operator function ${B}(t)$ in $\|\cdot\|_{1D}$ norm
the following equality holds:
$$
\gamma({B}(t))_{1D}= \gamma(D{B}(t)D^{-1})_{1}.
$$
Denote the elements of the matrix $D B(t)D^{-1}$ by $b_{ij}^*(t)$ i.e.
$D B(t)D^{-1}= \left(b_{ij}^*(t)\right)_{i,j=1}^{\infty}$.
Assume that
\begin{equation}
b_{ij}^*(t) \ge 0, \  i \neq j, \ t \ge 0.   \label{posit01}
\end{equation}
Put
\begin{equation}
\alpha_i\left(t\right)= \sum_{j=0}^\infty b_{ji}^*(t), \quad
\chi_i\left(t\right)= -\sum_{j=0}^\infty |b_{ji}^*(t)|, \ i \ge 1, \label{posit02}
\end{equation}
\noindent and let $\alpha(t)$
and $\beta(t)$ denote the
least lower and the least upper bound of the
sequence of functions $\{ \alpha_i(t), \ i \ge 1\}$
and $\chi$ denote the least upper bound of $\{ \chi_i(t), \ i \ge 1\}$
i.e.
\begin{equation}
\alpha\left(t\right)=\inf_{i \ge 1}\alpha_i\left(t\right), \ \beta
\left(t\right)=\sup_{i \ge 1}\alpha_i\left(t\right), \label{gl42.03}
\end{equation}
\noindent
\begin{equation}
\chi\left(t\right)=\sup_{i \ge 1}\chi_i\left(t\right).
\label{gl42.03''}
\end{equation}

Then the logarithmic norms of $B(t)$ and $(-B(t))$ are equal to
\begin{eqnarray}
\gamma \left(B\left(t\right)\right)_{1D} = \sup_{i} \alpha_i (t) = -
\alpha\left(t\right), \nonumber \label{gl42.10'}
\
\
\gamma \left(-B\left(t\right)\right)_{1D}
=   \sup \chi_i\left(t\right) =  \chi\left(t\right). \nonumber \label{gl42.10''}
\end{eqnarray}
\noindent If now one defines ${\bf v}\left(t\right) = D\left({\bf p^*}\left(t\right) -
{\bf p^{**}}\left(t\right)\right)$, then the following equation holds
\begin{equation}
\frac{d }{dt}{\bf v}(t) = DB\left(t\right)D^{-1} {\bf v}(t),
\label{gl42.047}
\end{equation}
\noindent
Notice that due to (\ref{posit01}), the inequality ${\bf v} \left(s\right) \ge {\bf 0}$ implies that ${\bf v}
\left(t\right) \ge {\bf 0}$ for any $t \ge s$. Hence
\begin{equation}
\frac{d} {dt} \sum_{i=1}^\infty v_i(t) \ge - \beta \left(t\right) \sum_{i=1}^\infty v_i(t),  \label{gl42.048}
\end{equation}
\noindent and one can obtain establish the following theorem.

\bigskip

\begin{theorem}\hspace{-0.2cm}
Let there exist an increasing sequence  $\{d_j, \ j \ge 1\}$ of positive
numbers with $d_1=1$, such that (\ref{posit01}) holds, and
$\alpha(t)$ defined by \eqref{gl42.03} satisfies
\begin{equation}
\int_0^{\infty} \alpha(t)\, dt = + \infty.
\end{equation}
Then the Markov chain  $\{ X(t), \ t\geq 0\}$ is weakly ergodic and the following bounds hold:
\begin{equation}
e^{-\int_s^t
{\chi\left(u\right)du}}\|\vp^*\left(s\right)-\vp^{**}\left(s\right)\|_{1D}
\le \|\vp^*\left(t\right)-\vp^{**}\left(t\right)\|_{1D} \le
e^{-\int_s^t
{\alpha\left(u\right)du}}\|\vp^*\left(s\right)-\vp^{**}\left(s\right)\|_{1D},
\label{t001}
\end{equation}
\begin{equation}
\|\vp^*(t)-\vp^{**}(t)\|  \le 4 e^{-\int_s^t
{\alpha(u)du}}\|\vz^*(s)-\vz^{**}(s)\|_{1D}, \label{t002}
\end{equation}
\begin{equation}
\|\vp^*(t)-\vp^{**}(t)\|_{1E} \le \frac{2}{W}e^{-\int_s^t
{\alpha(u)du}}\|\vz^*(s)-\vz^{**}(s)\|_{1D}, \label{t003}
\end{equation}
\noindent for any initial conditions $s \ge 0$, ${\vp^*}(s)$,
${\vp^{**}}(s)$ and any  $t \ge s$.

\smallskip
If in addition $D\left({\bf p^*}\left(s\right) - {\bf
p^{**}}\left(s\right)\right) \ge {\bf 0}$, then
\begin{equation}
\|\vp^*\left(t\right)-\vp^{**}\left(t\right)\|_{1D} \ge e^{-\int_s^t
{\beta\left(u\right)du}}\|\vp^*\left(s\right)-\vp^{**}\left(s\right)\|_{1D},
\label{t004}
\end{equation}
for any $0 \le s \le t$.
\end{theorem}

\noindent One can also obtain the corresponding
perturbation bounds. For the first results in this direction
see \cite{Kartashov1985,Kartashov1986,Kartashov1996,Zeifman1985},
for the stronger results see \cite{Mitrophanov2003} and
for the general approach see \cite{Zeifman2014s}.
The respective uniform in time truncation bounds can be obtained via techniques
proposed in \cite{Zeifman2014i,Zeifman2016t}.

If the Markov chain is  homogeneous, then all
elements $b_{ij}^*(t)$ of the matrix $D B(t)D^{-1}$ do not dependent on $t$
i.e. the quantities in \eqref{gl42.03} are constants.
Thus instead of general bounds given by Theorem 1,
one can specify then and obtain the following theorem.

\begin{theorem}\hspace{-0.2cm}
Let there exist an increasing sequence  $\{d_j, \ j \ge 1\}$ of positive
numbers with $d_1=1$, such that (\ref{posit01}) holds, and
$\alpha(t)=\alpha$ defined by \eqref{gl42.03} is positive i.e.
$\alpha > 0$.
%$\{d_j\}$ of positive
%numbers  such that  $d_1=1$, (\ref{posit01}) hold, and in addition,
%\begin{equation}
%\alpha > 0.
%\end{equation}
Then the Markov chain  $\{ X(t), \ t\geq 0\}$ is strongly ergodic and the following bounds hold:
\begin{equation}
e^{- \chi t}\|\vp^*\left(0\right)-\vp^{**}\left(0\right)\|_{1D} \le
\|\vp^*\left(t\right)-\vp^{**}\left(t\right)\|_{1D} \le e^{- \alpha
t}\|\vp^*\left(0\right)-\vp^{**}\left(0\right)\|_{1D},
\end{equation}
\begin{equation}
\|\vp^*(t)-\vp^{**}(t)\|  \le 4 e^{- \alpha
t}\|\vz^*(0)-\vz^{**}(0)\|_{1D},
\end{equation}
\begin{equation}
\|\vp^*(t)-\vp^{**}(t)\|_{1E} \le \frac{2}{W}e^{-\alpha
t}\|\vz^*(0)-\vz^{**}(0)\|_{1D},
\end{equation}
\noindent for any initial conditions $s \ge 0$, ${\vp^*}(0)$,
${\vp^{**}}(0)$ and any  $t \ge 0$.
\smallskip
If in addition $D\left({\bf p^*}\left(0\right) - {\bf
p^{**}}\left(0\right)\right) \ge {\bf 0}$, then
\begin{equation}
\|\vp^*\left(t\right)-\vp^{**}\left(t\right)\|_{1D} \ge e^{- \beta
t}\|\vp^*\left(0\right)-\vp^{**}\left(0\right)\|_{1D},
\end{equation}
for any $ t \ge 0$.

For the decay parameter $\alpha^*$ defined as
$$
\lim_{t \rightarrow \infty} (p_{ij}(t)-\pi_j)= O(e^{-\alpha^* t}),
$$
where $\{ \pi_j, \ j \ge 0\}$ are the stationary probabilities
of the chain, it holds that $\alpha^* \ge \alpha$.
\end{theorem}

\noindent Notice that some additional results related to \textit{Theorem} 2 can
also be found in
\cite{Doorn2010,Granovsky2000,Zeifman1991}.
If one assumes that the intensities $q_{ij}(t)$ are $1-$periodic in $t$ i.e.
$q_{ij}(t)$ are periodic functions and the length of the period is equal to one,
then the Markov chain  $\{ X(t), \ t\geq 0\}$ has the limiting $1-$periodic
limiting regime. Under the assumptions of Theorem 1 the Markov chain  $\{ X(t), \ t\geq 0\}$
 is exponentially weakly ergodic. The detailed discussion of this results is
 given in \cite{Zeifman2006}.

Till the end of this section one presents
a bit more detailed analysis of two special cases:
homogeneous case and the case with periodic
intensities. Firstly note that in the both cases there
exist positive $R$ and $a$
such that
\begin{equation}
e^{-\int_s^{t} \alpha(u)\, du} \le R e^{-a(t-s)} \label{311}
\end{equation}
\noindent for any $0 \le s \le t$. Hence  the Markov chain  $\{ X(t), \ t\geq 0\}$
is {\it exponentially} weakly ergodic.
Indeed, if the Markov chain  $\{ X(t), \ t\geq 0\}$ is homogeneous,
then one may put $R=1$, $a = \alpha$ given by \eqref{gl42.03}.
If all the intensity functions $q_{ij}(t)$  are $1-$periodic in $t$,
then one may put
$$a = \int_0^1 \alpha(t)\, dt, \quad R = e^K,
\quad K = \sup_{|t-s| \le 1} \int_s^t\alpha(u)\, du.$$
By doing so, for any solution of (\ref{2.06}) the following bound holds:
\begin{eqnarray} \label{317}
\|{\bf z}(t)\|_{1D} \le  \nonumber \\ \|V(t)\|_{1D}\|{\bf
z}(0)\|_{1D} + \int_0^t
\|V(t,\tau)\|_{1D} \|{\bf f} (\tau)\|_{1D}\, d \tau  \le  \\
R e^{-a t}\|{\bf z}(0)\|_{1D} + \frac{F R}{a},\nonumber
\end{eqnarray}
\noindent where $F$ is such that $\|{\bf f} (t)\|_{1D} \le F$ for
almost all $t \in [0,1]$.
Hence one has the upper bound for the limit
\begin{equation}
\limsup_{t \to \infty}  \|{\bf z}(t)\|_{1D} \le  \frac{F R}{a},
\label{318}
\end{equation}
\noindent for any initial condition and
\begin{eqnarray} \|{\bf p}(0) - {\bf e}_0\|_{1D} =
\|{\bf p}(0)\|_{1D} = \|{\bf z}(0)\|_{1D} \le \limsup_{t \to
\infty} \|{\bf z}(t)\|_{1D}, \label{343}
\end{eqnarray}
where ${\bf e}_i$ denotes the unit vector of zeros with 1 in the $i$-th place.
If the initial distribution is $\vp^{**}(0) = {\bf e}_0$ then  $\vz^{**}(0) = {\bf 0}$,
$\vz(t) \ge  {0}$ for any $\vp^{*}(0)$ and any $t \ge 0$.
Therefore
\begin{eqnarray}
\|\vz(t)\|_{1D}= d_1p_1+(d_1+d_2)p_2+(d_1+d_2+d_3)p_3+...=
\nonumber
\\ = d_1
p_1+\frac{d_1+d_2}{2} 2p_2+\frac{d_1+d_2+d_3}{3} 3p_3+...\ge
\inf_k\frac{d_1+...+d_k}{k} \|{\bf z}(t)\|_{1E},\nonumber
\end{eqnarray}
\noindent and one can use $W^*=\inf_k\frac{d_1+...+d_k}{k}$ instead
of  $W=\inf_k \frac {d_k}{k}$, given by \eqref{2002}
in all the bounds on the rate of
convergence.
Finally, for the considered two special cases one has the following two corollaries.

\smallskip

\begin{corollary}
Let  $\{ X(t), \ t\geq 0\}$ be a homogeneous Markov chain  and let there exist an
increasing sequence $\{d_j, \ j \ge 1\}$ of positive numbers with
$d_1=1$ such that (\ref{posit01}) holds and in addition $\alpha > 0$.
Then the Markov chain
$\{ X(t), \ t\geq 0\}$ is exponentially ergodic
and the following bounds hold:
\begin{equation}
\|{\bf \pi} -{\bf p}(t,0)\|  \le \frac{4F}{\alpha} e^{-\alpha t},
\label{cor001}
\end{equation}
\begin{equation}
|\phi -E(t,0)| \le \frac{F}{\alpha W^*} e^{-\alpha t},
\label{cor002}
\end{equation}
\noindent where
${\pi} = (\pi_0, \pi_1, \dots)^T$ denotes the vector of stationary probabilities
of the chain and $\phi= \sum_{j=0}^\infty j \pi_j$ and ${\bf p}(0,0)={\bf e}_0$.
\end{corollary}

\begin{corollary}
Assume all intensity function of
the Markov chain $\{ X(t), \ t\geq 0\}$ be $1-$periodic in $t$.
Let there exist an increasing sequence $\{d_j, \ j \ge 1\}$ of positive numbers
with $d_1=1$ such that (\ref{posit01}) holds and in addition $ \int_0^{1}
\alpha(t)\, dt = a > 0$. Then the Markov chain $\{ X(t), \ t\geq 0\}$
is exponentially weakly ergodic
and the following bounds hold:
\begin{equation}
\|{\bf \pi}(t) -{\bf p}(t,0)\|  \le \frac{4FR}{a} e^{-a t},
\label{cor003}
\end{equation}
\begin{equation}
|\phi(t) -E(t,0)| \le \frac{FR}{aW^*} e^{-a t}, \label{cor004}
\end{equation}
\noindent where
${\pi} (t)= (\pi_0(t), \pi_1(t), \dots)^T$ denotes the vector of
limiting probabilities
of the chain
and $\phi(t)= \sum_{j=0}^\infty j \pi_j(t)$ and ${\bf p}(0,0)={\bf
e}_0$.
\end{corollary}

\noindent If the state space of the Markov chain is finite
there exist a number of special results (one can refer to \cite{Doorn2010,Granovsky2000,Zeifman2013m}).

\section{Convergence bounds}

In order to apply the results of the \textit{Theorem~1} and
\textit{Theorem~2} and to obtain the convergence bounds
for the system from either ${\rm I}^{st}-{\rm IV}^{th}$ class,
one has to know the exact expressions
for the functions $\alpha_i(t)$ and $\chi_i\left(t\right)$,
given by \eqref{gl42.03} and \eqref{gl42.03''}.
In this section one provides
the expressions for
$\alpha_i(t)$ and $\chi_i\left(t\right)$
for single server systems from
classes ${\rm I}^{st}-{\rm IV}^{th}$.
Then one shows how these expression
change, when one considers to multiple server case.

\subsection{Inhomogeneous $M/M/S$ queueing system with batch arrivals and state-dependent service intensity}

Consider the queueing system $M/M/1$ queueing system
with time-dependent arrival and service intensities.
Let $\lambda_k(t)$ be the arrival intensity of the batch, containing
$k$ customers, at instant $t$  and $\mu_k(t)$ be the service
intensity at instant $t$ if the total number of
customers in the system is equal to $k$.
Then the transposed intensity matrix has the
form
\begin{equation}
A(t) = \begin{pmatrix} a_{00}\left(t\right) & \mu_1\left(t\right) &
0 & 0 & 0 & \cdots  \cr \la_1\left(t\right) & a_{11}\left(t\right) &
\mu_2\left(t\right) & 0 & 0& \cdots \cr \la_2\left(t\right) &
\la_1\left(t\right) & a_{22}\left(t\right) & \mu_3\left(t\right) &
0& \cdots \cr \la_3 \left(t\right) & \la_2 \left(t\right) & \la_1
\left(t\right) & a_{33}\left(t\right) & \mu_4\left(t\right) & \cdots
\cr \la_4\left(t\right) & \la_3\left(t\right) & \la_2\left(t\right)
& \la_1\left(t\right) & a_{44}\left(t\right) & \cdots  \cr \cdots &
\vdots& \vdots& \vdots& \vdots& \ddots&   \cr
\end{pmatrix},
\label{2.001}
\end{equation}
\noindent where diagonal elements of $A(t)$ are such that all column
sums are equal to zero for any $t \ge 0$. As the assumption
(\ref{posit01}) is fulfilled, it holds that
\begin{equation}
\alpha_j(t)=
\mu_j\left(t\right)-\frac{d_{j-1}}{d_j}\mu_{j-1}\left(t\right)+\sum_{i=1}^{\infty} \left(1-\frac{d_{i+j}}{d_j}\right)\la_i\left(t\right),
\label{3.02}
\end{equation}
\noindent and
\begin{equation}
\chi_j\left(t\right)= \mu_j\left(t\right)+\frac{d_{j-1}}{d_j}\mu_{j-1}\left(t\right)+\sum_{i=1}^{\infty} \left(1+\frac{d_{i+j}}{d_j}\right)\la_i\left(t\right).
\label{3.02''}
\end{equation}
\noindent
Therefore, the \textit{Theorem 1} and \textit{Theorem 2} hold for these $\alpha_i(t)$ and
$\chi_i(t)$.

Consider now the queueing system $M/M/S$ queueing system with $S>1$ servers,
 time-dependent arrival and service intensities. Customers arrive at the system
 in batches of size not greater than $S$.
 Assume that the arrival intensity of batch, containing $k$ customers, at instant $t$ is equal to
 $\lambda_{k}(t)=\frac{1}{S k}\lambda (t)$ if $1 \le k \le S$
 and $\lambda_{k}(t)=0$ if $k >S$.
 Denote the service intensity at instant $t$ by $\mu_{k}(t)$
 and assume that  $\mu_{k}(t)=\min{(k,S)}\mu (t)$.
{For the assumed values of $\lambda_{k}(t)$ and $\mu_{k}(t)$,
and the expressions for $\alpha_i(t)$ and $\chi_i(t)$ given above,
the \textit{Theorem 1} and \textit{Theorem 2} hold.
}

\subsection{Inhomogeneous $M/M/S$ queueing system with batch service and state-dependent arrival intensity}

Consider the queueing system $M/M/1$ queueing system
with time-dependent arrival and service intensities.
But now
let $\lambda_k(t)$ be the arrival intensity of
$k$ customers at instant $t$
if the total number of
customers in the system is equal to $k$
and $\mu_k(t)$ be the service
intensity at instant $t$ of a group of $k$ customers.
Then the transposed intensity matrix has the form
\begin{equation}
A(t)=\left(
\begin{array}{ccccccc}
a_{00}(t) & \mu_1(t)  & \mu_2(t)   & \mu_3(t)  & \mu_4(t) & \mu_5(t) & \cdots \\
\lambda_0(t)   & a_{11}(t)  & \mu_1(t)  & \mu_2(t)   & \mu_3(t)& \mu_4(t)  & \cdots  \\
0  & \lambda_1 (t)    & a_{22}(t)& \mu_1(t)  & \mu_2(t) & \mu_3(t)   &  \cdots \\
0& 0  & \lambda_2 (t)    & a_{33}(t)& \mu_1(t)  & \mu_2(t)    &  \cdots \\
\vdots& \vdots  & \vdots    & \vdots & \vdots  & \vdots    &  \ddots \\
\end{array}
\right),
\end{equation}
\noindent where
$a_{ii}\left(t\right)=-\sum_{k=1}^{i}\mu_k\left(t\right) -
\la_{i}\left(t\right)$.
Then \textit{Theorem 1} and \textit{Theorem 2} hold for
\begin{equation}
\alpha_i\left(t\right) = \mu_i\left(t\right) - \sum_{k=1}^{i-1}\left(\mu_{i-k}\left(t\right)-\mu_i\left(t\right)\right)\frac{d_k}{d_i}+\la_{i-1}\left(t\right)-\frac{d_{i+1}}{d_i}\la_{i}\left(t\right),
\end{equation}
\noindent
\begin{equation}
\chi_i\left(t\right)= \mu_i\left(t\right) + \sum_{k=1}^{i-1}\left(\mu_{i-k}\left(t\right)-\mu_i\left(t\right)\right)\frac{d_k}{d_i}+\la_{i-1}\left(t\right)+\frac{d_{i+1}}{d_i}\la_{i}\left(t\right).
\end{equation}

Consider again the queueing system $M/M/S$ queueing system with $S>1$ servers,
 time-dependent arrival and service intensities. The customers are served
 in batches of size not greater than $S$.
 Assume that the service intensity is
$\mu_{k}(t)=\frac{1}{k}\mu (t)$ if  $1 \le k \le S$
and $\mu_{k}(t)=0$ if $k>S$.
Denote the arrival intensity by $\lambda (t)$.
{For the assumed values of $\lambda_{k}(t)$ and $\mu_{k}(t)$,
and the expressions for $\alpha_i(t)$ and $\chi_i(t)$ given above,
the \textit{Theorem 1} and \textit{Theorem 2} hold.
}

\subsection{Inhomogeneous $M/M/S$ queueing system with batch arrivals and service}

Consider the queueing system $M/M/1$ queueing system
with time-dependent arrival and service intensities.
Customer arrive and are served in batches.
Let $\lambda_k(t)$ and $\mu_k(t)$
be the arrival and service intensity of a group of $k$ customers.
This queueing system has been extensively studied
with respect to the rate of convergence, truncation and perturbation bounds
in \cite{Satin2013,Zeifman2014a}.
The transposed intensity matrix in this case
has the following form:
\begin{equation}
A\left(t\right)=\left(
\begin{array}{ccccccc}
a_{00}\left(t\right) & \mu_1\left(t\right)  & \mu_2\left(t\right)   & \mu_3\left(t\right)  & \cdots  \\
\la_1\left(t\right)   & a_{11}\left(t\right)  & \mu_1\left(t\right)  & \mu_2\left(t\right)    & \cdots  \\
\la_2\left(t\right)  & \la_1\left(t\right)    & a_{22}\left(t\right)& \mu_1\left(t\right)    &  \cdots  \\
\la_3\left(t\right)  & \la_2\left(t\right)  & \la_1\left(t\right)    & a_{33}\left(t\right)&  \cdots    \\
&  \vdots  &  \vdots  &  \vdots  &  \ddots
\end{array}
\right),
\end{equation}
\noindent and
$a_{ii}\left(t\right)=-\sum_{k=1}^{i}\mu_k\left(t\right) -
\sum_{k=1}^{\infty} \la_{k}\left(t\right)$.
Therefore, \textit{Theorem 1} and \textit{Theorem 2} hold for
\begin{equation}
\alpha_i\left(t\right)=- a_{ii}\left(t\right) - \sum_{k=1}^{i-1}\left(\mu_{i-k}\left(t\right)-\mu_i\left(t\right)\right)
\frac{d_k}{d_i}-\sum_{k \ge 1} \la_k\left(t\right) \frac{d_{k+i}}{d_i},
\end{equation}
\noindent and
\begin{equation}
\chi_i\left(t\right)= -a_{ii}\left(t\right) + \sum_{k=1}^{i-1}\left(\mu_{i-k}\left(t\right)-
\mu_i\left(t\right)\right)\frac{d_k}{d_i}+ \sum_{k \ge 1}\la_k\left(t\right)\frac{d_{k+i}}{d_i}.
\end{equation}

Consider again the queueing system $M/M/S$ queueing system with $S>1$ servers,
 time-dependent arrival and service intensities.
 Assume the arrivals and service appear in batches of size
 not greater than $S$.
 Assume that the arrival intensity of $k$ customers at instant $t$ is equal to
 $\lambda_{k}(t)=\frac{1}{S k}\lambda (t)$ if $1 \le k \le S$
 and $\lambda_{k}(t)=0$ if $k >S$;
 the service intensity is assume to be equal to
$\mu_{k}(t)=\frac{1}{k}\mu (t)$ if  $1 \le k \le S$
and $\mu_{k}(t)=0$ if $k>S$.
{For the assumed values of $\lambda_{k}(t)$ and $\mu_{k}(t)$,
and the expressions for $\alpha_i(t)$ and $\chi_i(t)$ given above,
the \textit{Theorem 1} and \textit{Theorem 2} hold.
}

\subsection{Inhomogeneous $M/M/S$  queueing system with state-dependent arrival and service intensities}

If in the queueing system $M/M/1$ queueing system
the arrival $\lambda_n(t)$  and service intensities $\mu_n(t)$
are time-dependent
and also depend on the total number of customers $n$
in the system
then the queue-length process for a general Markovian queue
with intensities  $\lambda_n(t)$ is
the inhomogeneous birth-death process
with birth and death intensities equal to $\lambda_n(t)$  and $\mu_n(t)$ .
Thus  \textit{Theorem 1} and \textit{Theorem 2} hold with
\begin{equation}
\alpha_j(t)=
\mu_{j}\left(t\right)-\frac{d_{j-1}}{d_j}\mu_{j-1}\left(t\right)+\la_{j-1}\left(t\right)-\frac{d_{j+1}}{d_j}\la_{j}\left(t\right),
\end{equation}
\noindent and
\begin{equation}
\chi_j(t)= \mu_{j}\left(t\right)+\frac{d_{j-1}}{d_j}\mu_{j-1}\left(t\right)+\la_{j-1}\left(t\right)+\frac{d_{j+1}}{d_j}\la_{j}\left(t\right).
\end{equation}

Consider again the ordinary queueing system $M/M/S$ queueing system with $S>1$ servers,
 time-dependent arrival and service intensities $\lambda(t)$ and $\mu_n(t) = \min(n,S)\mu(t)$
 correspondingly.
{For the assumed values of $\lambda_{k}(t)$ and $\mu_{k}(t)$,
and the expressions for $\alpha_i(t)$ and $\chi_i(t)$ given above,
the \textit{Theorem 1} and \textit{Theorem 2} hold.
}

\section{Numerical examples}

The purpose of the numerical section is two-fold.
Firstly one
demonstrates that the convergence bounds obtained
in the previous section can indeed be computed.
Having fixed the arrival and service intensities in the inhomogeneous $M/M/S$ queueing system
with state-independent arrival and service intensities
($I^{st}$ class),
one specifies the sequence $\{ d_i, \ i \ge 1\}$ and
provides corresponding bounds
using \textit{Corollary 2}.
Secondly one shows that the approach
proposed in this paper can be used to
compute approximations for the
limiting characteristics of the systems
with a given approximation error.
The characteristics under consideration are the
limiting idle probability and the limiting mean number of customers
in the system.

The system considered are:

(i) inhomogeneous $M/M/100$ queueing system with state-independent arrival and service intensities;

(ii) inhomogeneous $M/M/100$ queueing system with state-independent batch arrivals;

(iii) inhomogeneous $M/M/100$ queueing system with state-independent batch service.

(iv) inhomogeneous $M/M/100$ queueing system with state-independent batch arrivals and batch service.

All transition intensities are assumed to be periodic functions of time.
Customer in all three systems are served in FCFS order.
The inhomogeneous $M/M/S$ consists of single infinite capacity queue
and $100$ servers. Arrivals happen according
to the inhomogeneous Poisson process with  the intensity $\lambda^*(t;i)$
equal to
$$
\lambda^*(t;i)= i(1+\sin 2\pi t), \ i>0, \ t >0.
$$
\noindent Whenever the server
becomes free, a customer from the queue (if there is any) enters server and
get served according to exponential distribution
with the intensity
$$
\mu^*(t)= 3+\cos 2\pi t.
$$

In the inhomogeneous $M/M/100$ queueing system with batch arrivals
it is assumed that customers arrive
in batches in accordance with a inhomogeneous Poisson process
of intensity $\lambda^*(t;i) / \sum_{i=1}^S (Si)^{-1}$.
The size of the arriving group is a random variable
with discrete probability distribution $(S k)^{-1} / \sum_{i=1}^S (Si)^{-1}$, $1 \le k\le S$.
The sizes and interarrival times of successive arriving groups are
stochastically independent. Thus the total arrival intensity
is $\lambda_{k}(t,i)=\frac{1}{S k}\lambda^*(t,i)$ if $1 \le k\le S$
and $\lambda_{k}(t,i)=0$ if $k >S$.
Whenever the server
becomes free, a customer from the queue (if there is any) enters server and
get served according to exponential distribution
with the intensity $\mu_{k}(t)=\min{(k,S)}\mu^*(t)$.

In the inhomogeneous $M/M/100$ queueing system with batch service
customers arrive
in accordance with a inhomogeneous Poisson process
of the same intensity $\lambda^*(t,i)$.
But the service happens in batches of size not
greater than $S$ and the service times
are exponentially distribution with the service intensity
equal to $\mu_{k}(t)=\frac{1}{k}\mu (t)$ if  $1 \le k \le S$
and $\mu_{k}(t)=0$ if $k>S$.

Finally, in the inhomogeneous $M/M/100$ queueing system with batch arrivals and
batch service, the
total arrival intensity
is $\lambda_{k}(t,i)=\frac{1}{S k}\lambda^*(t,i)$ if $1 \le k\le S$
and $\lambda_{k}(t,i)=0$ if $k >S$
and the service intensity is
$\mu_{k}(t)=\frac{1}{k}\mu (t)$ if  $1 \le k \le S$
and $\mu_{k}(t)=0$ if $k>S$.

Let us find the convergence bounds in case (i) i.e. for
the inhomogeneous $M/M/100$ queueing system with state-independent arrival and service intensities.
Let $i=50$ i.e. arrival intensity is equal to
$$
\lambda^*(t;50)= 50(1+\sin 2\pi t), \ t >0.
$$
Specify the sequence $\{ d_i, \ i \ge 1\}$
as follows:

-- $d_i=1$, if $1 \le i \le 100$;

-- $d_{101}=1.05d_{100}$, $d_{102}=1.1d_{101}$, $d_{103}=1.3d_{102}$, $d_{104}=1.6d_{103}$,
$d_{105}=2d_{104}$;

-- $d_{i}=2.3^{i-105} d_{105}$, if $i \ge 106$.

\noindent Such sequence $\{ d_i, \ i \ge 1\}$ guarantees that the
assumptions of \textit{Corollary~2} are fulfilled and one has
bounds on the rate of convergence to the limiting characteristics
given by (\ref{cor003}) and (\ref{cor004}) with $a=1.7$, $R=2$, $F=100$.

In order to approximate the limiting characteristics for
all three cases (i)--(iv),
one can apply \textit{Theorem~5} and \textit{Theorem~8} from \cite{Zeifman2014i}.
But firstly one has to specify the value $i$ in the arrival intensity $\lambda^*(t;i)$,
because as the arrival intensity (and thus load) grows the
bigger state space is needed.
Assume that $i\le 50$ i.e. the maximum arrival intensity allowed
is $\lambda^*(t;50)= 50(1+\sin 2\pi t)$.
Then one can compute the solution of the forward Kolmogorov system
(\ref{ur01}) for the truncated process on the state space $\{
0, 1, \dots, 155 \}$ on the interval $[0,t^*+1]$ with
the initial condition $X(0)=0$. Hence one finds the limiting
idle probability and limiting mean value on the interval $[t^*,t^*+1]$ with an
error less than $10^{-4}$, where $t^*=5$ or $t^*=7$.

Below one presents the plots of the limiting probability of the empty
queue $p_0(t)$ and the limiting mean $\varphi (t)$ number of customer in the system
for all each of the cases (i)--(iv).

%%%%%%% case (i)

\begin{center}
\begin{figure}[!ht]
\begin{center}
  \includegraphics[scale=0.35]{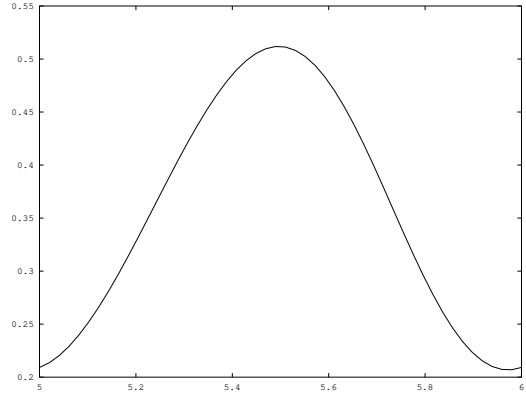}
\end{center}
\caption{Case (i), the arrival intensity is $\lambda^*(t;10)$. Approximation of the limiting
mean number $\varphi (t)$ of customers in the system for $t\in[5,6]$.}
\label{fig:41}       % Give a unique label
\end{figure}
\end{center}

\begin{center}
\begin{figure}[!ht]
\begin{center}
  \includegraphics[scale=0.35]{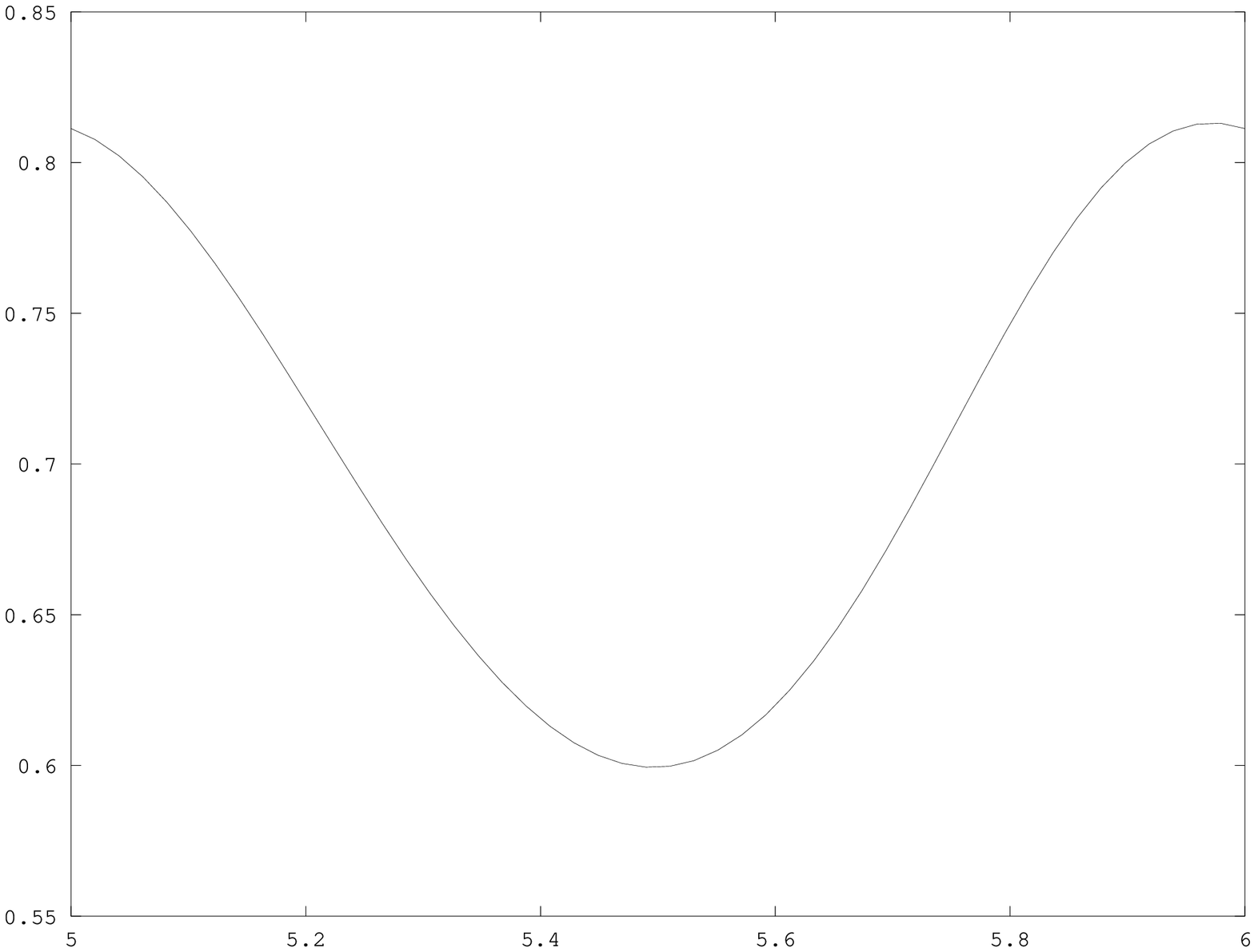}
\end{center}
\caption{Case (i), the arrival intensity is $\lambda^*(t;10)$. Approximation of the limiting probability $p_0(t)$
of the empty queue for $t\in[5,6]$.}
\label{fig:42}       % Give a unique label
\end{figure}
\end{center}

\begin{center}
\begin{figure}[!ht]
\begin{center}
  \includegraphics[scale=0.35]{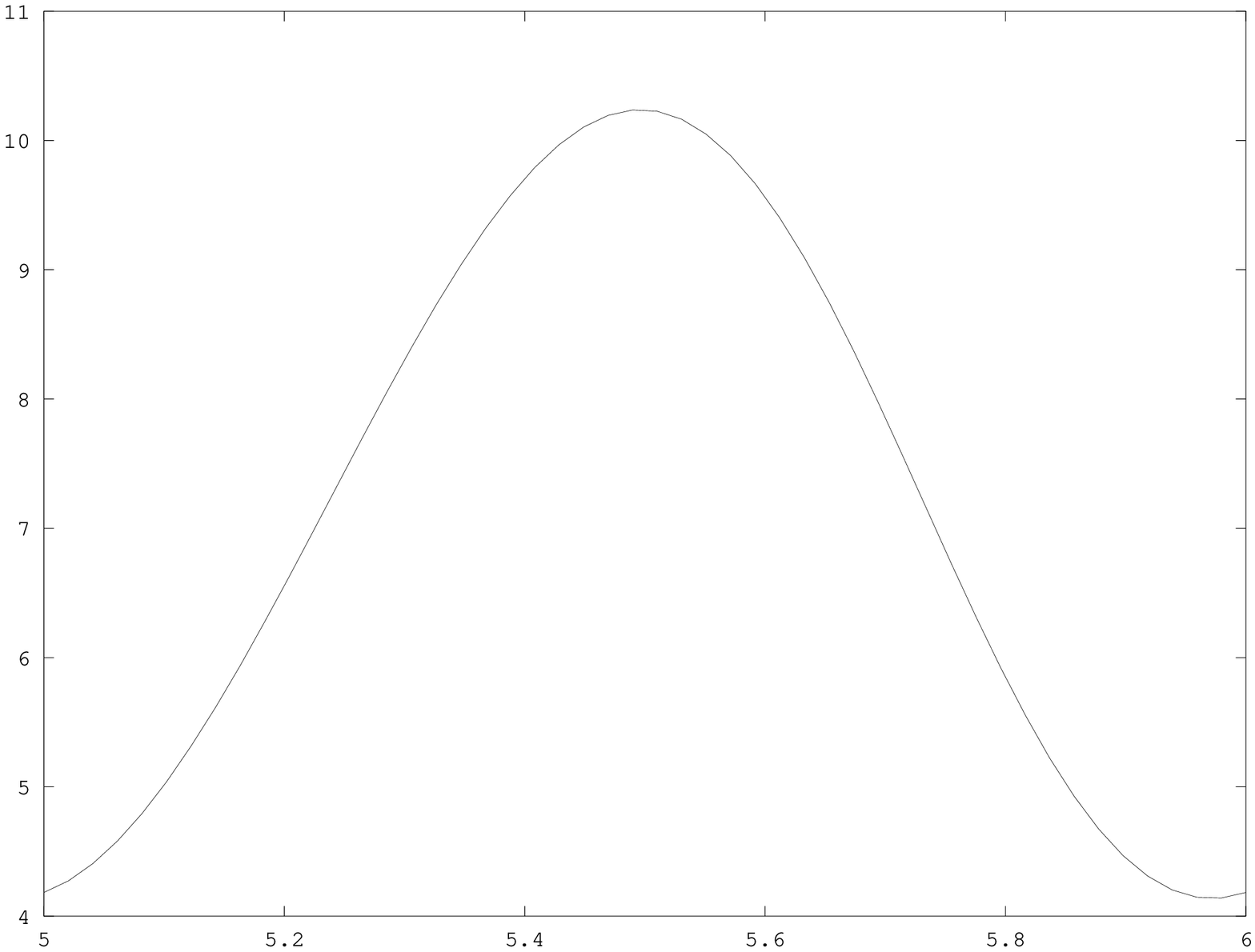}
\end{center}
\caption{Case (i), the arrival intensity is $\lambda^*(t;20)$. Approximation of the limiting
mean number $\varphi (t)$ of customers in the system for $t\in[5,6]$.}
\label{fig:43}       % Give a unique label
\end{figure}
\end{center}

\begin{center}
\begin{figure}[!ht]
\begin{center}
  \includegraphics[scale=0.35]{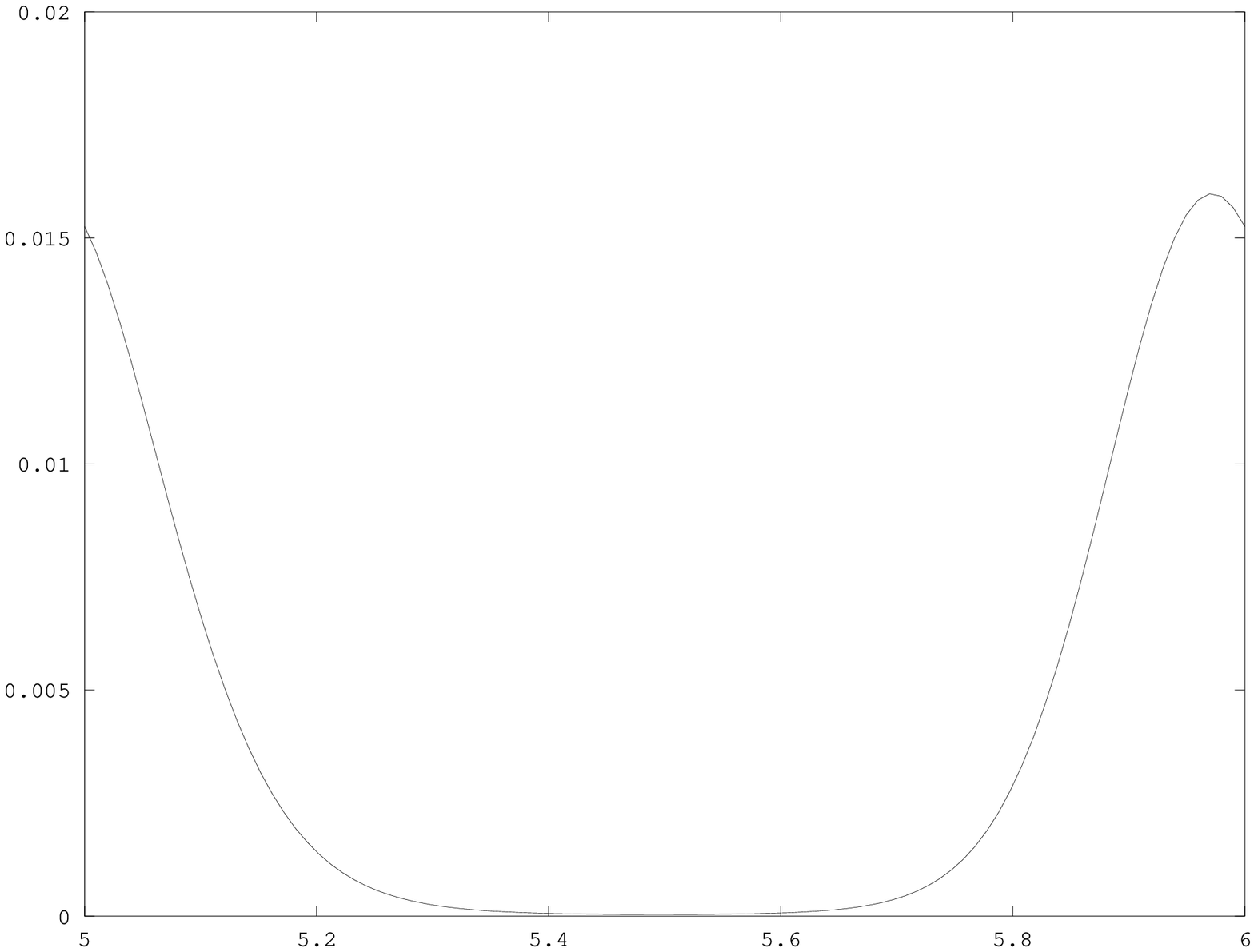}
\end{center}
\caption{Case (i), the arrival intensity is $\lambda^*(t;20)$. Approximation of the limiting probability  $p_0(t)$
of the empty queue for $t\in[5,6]$.}
\label{fig:44}       % Give a unique label
\end{figure}
\end{center}

\begin{center}
\begin{figure}[!ht]
\begin{center}
  \includegraphics[scale=0.35]{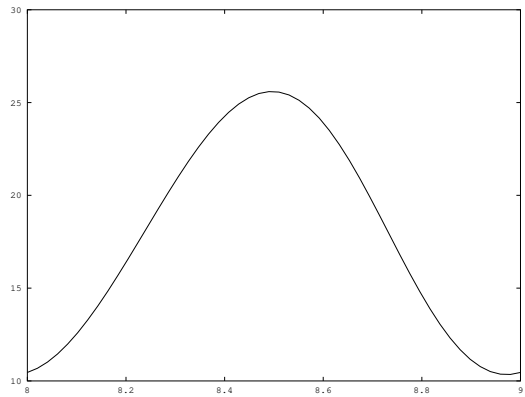}
\end{center}
\caption{Case (i), the arrival intensity is $\lambda^*(t;50)$. Approximation of the limiting
mean number $\varphi (t)$ of customers in the system for $t\in[5,6]$.}
\label{fig:45}       % Give a unique label
\end{figure}
\end{center}

\clearpage

\noindent Note from the figures above, that in case of high arrival intensity the
limiting probability $p_0(t)$  of the empty queue here equals to 0
most of the time.

%%%%%%%%%% case (ii)

\begin{center}
\begin{figure}[!ht]
\begin{center}
  \includegraphics[scale=0.35]{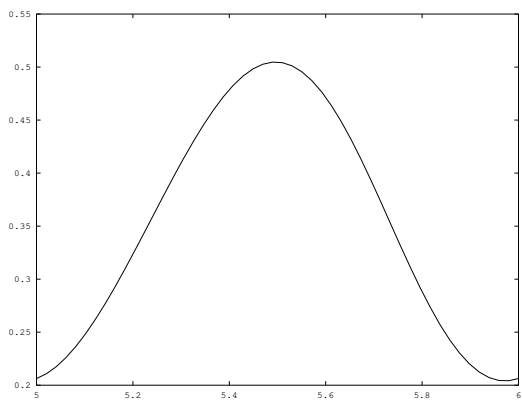}
\end{center}
\caption{Case (ii), the arrival intensity is $\lambda^*(t;10)$. Approximation of the limiting
mean number $\varphi (t)$ of customers in the system for $t\in[5,6]$.}
\label{fig:1}       % Give a unique label
\end{figure}
\end{center}

\begin{center}
\begin{figure}[!ht]
\begin{center}
  \includegraphics[scale=0.35]{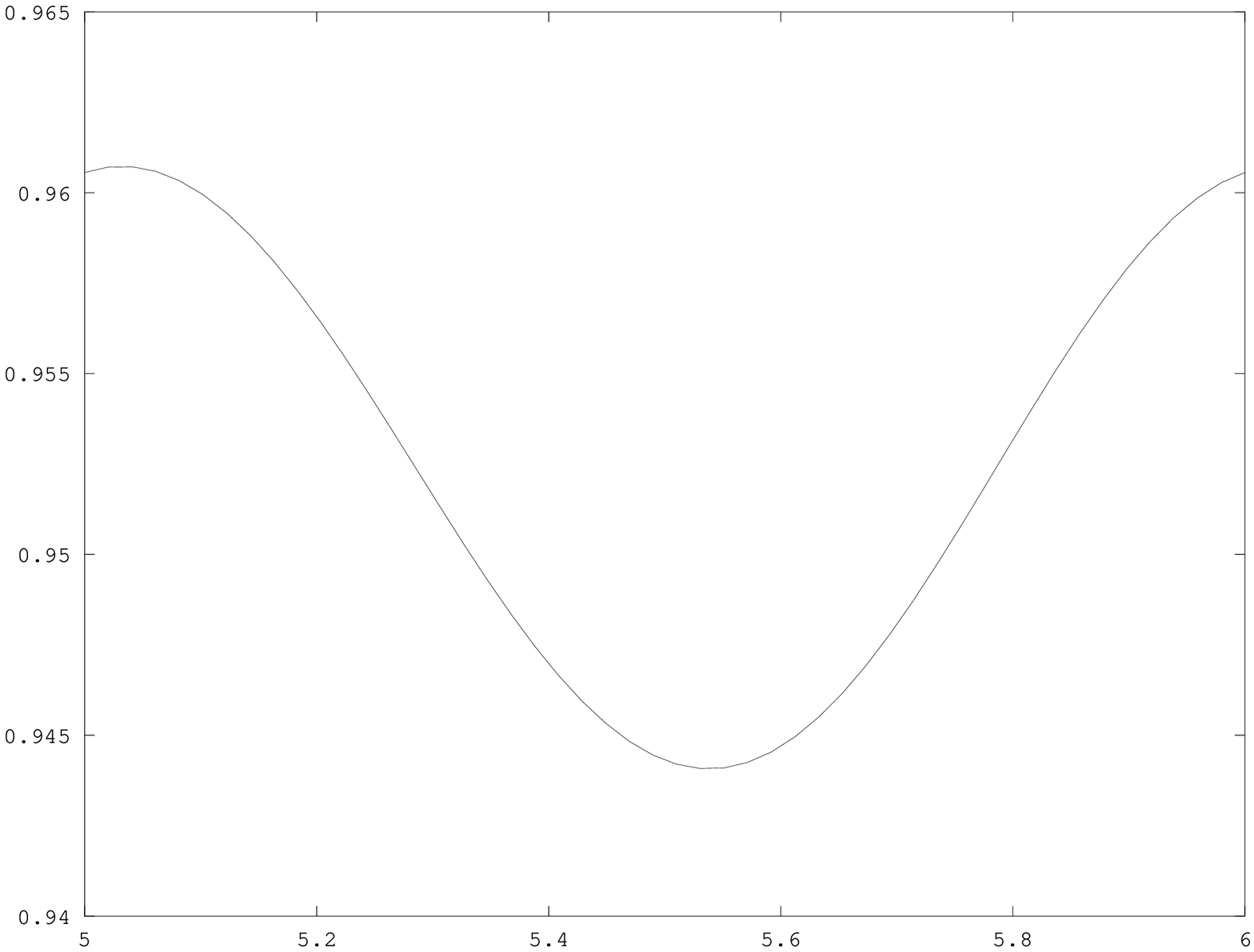}
\end{center}
\caption{Case (ii), the arrival intensity is $\lambda^*(t;10)$. Approximation of the limiting probability $p_0(t)$
of the empty queue for $t\in[5,6]$.}
\label{fig:2}       % Give a unique label
\end{figure}
\end{center}

\begin{center}
\begin{figure}[!ht]
\begin{center}
  \includegraphics[scale=0.35]{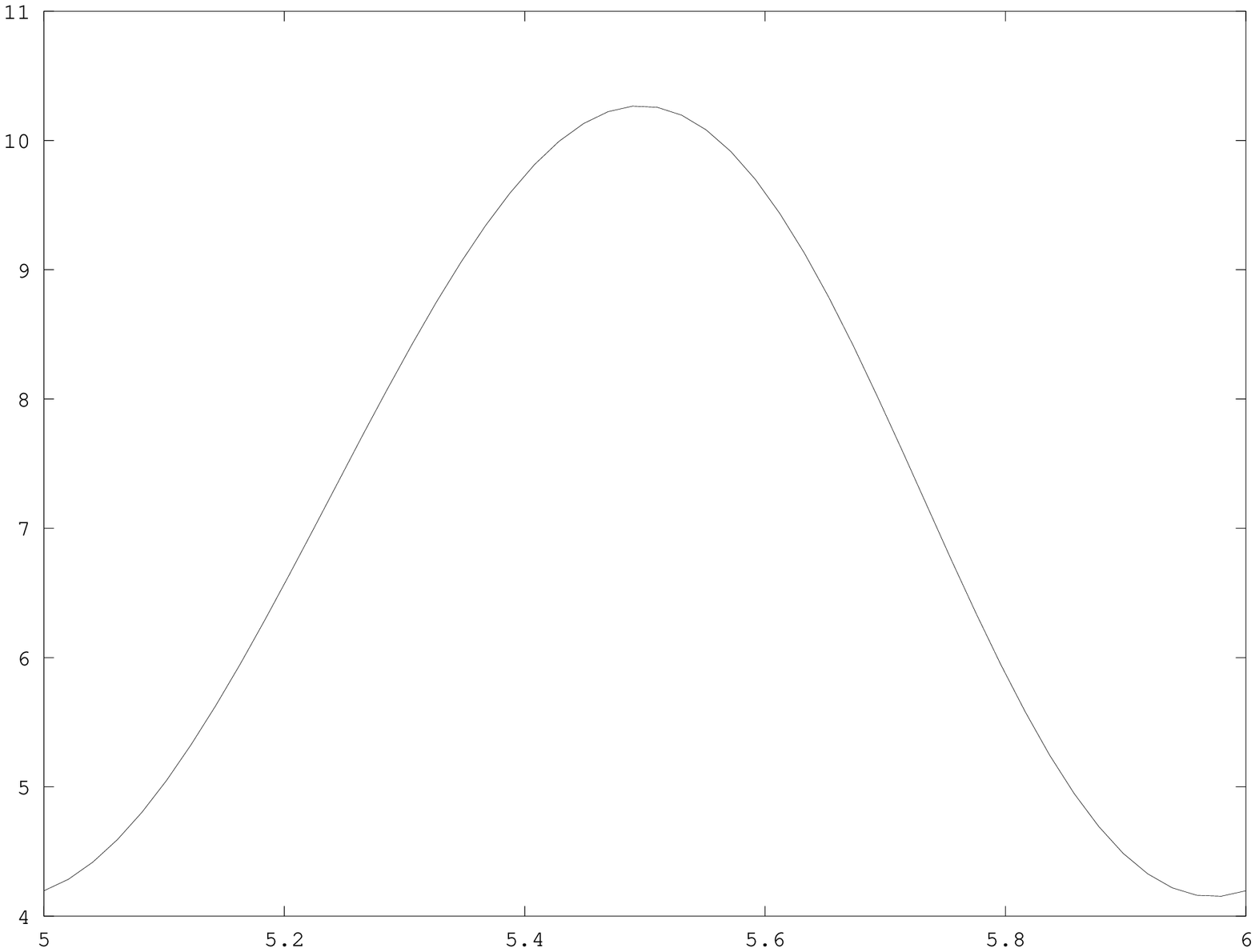}
\end{center}
\caption{Case (ii), the arrival intensity is $\lambda^*(t;20)$. Approximation of the limiting
mean number $\varphi (t)$ of customers in the system for $t\in[5,6]$.}
\label{fig:3}       % Give a unique label
\end{figure}
\end{center}

\begin{center}
\begin{figure}[!ht]
\begin{center}
  \includegraphics[scale=0.35]{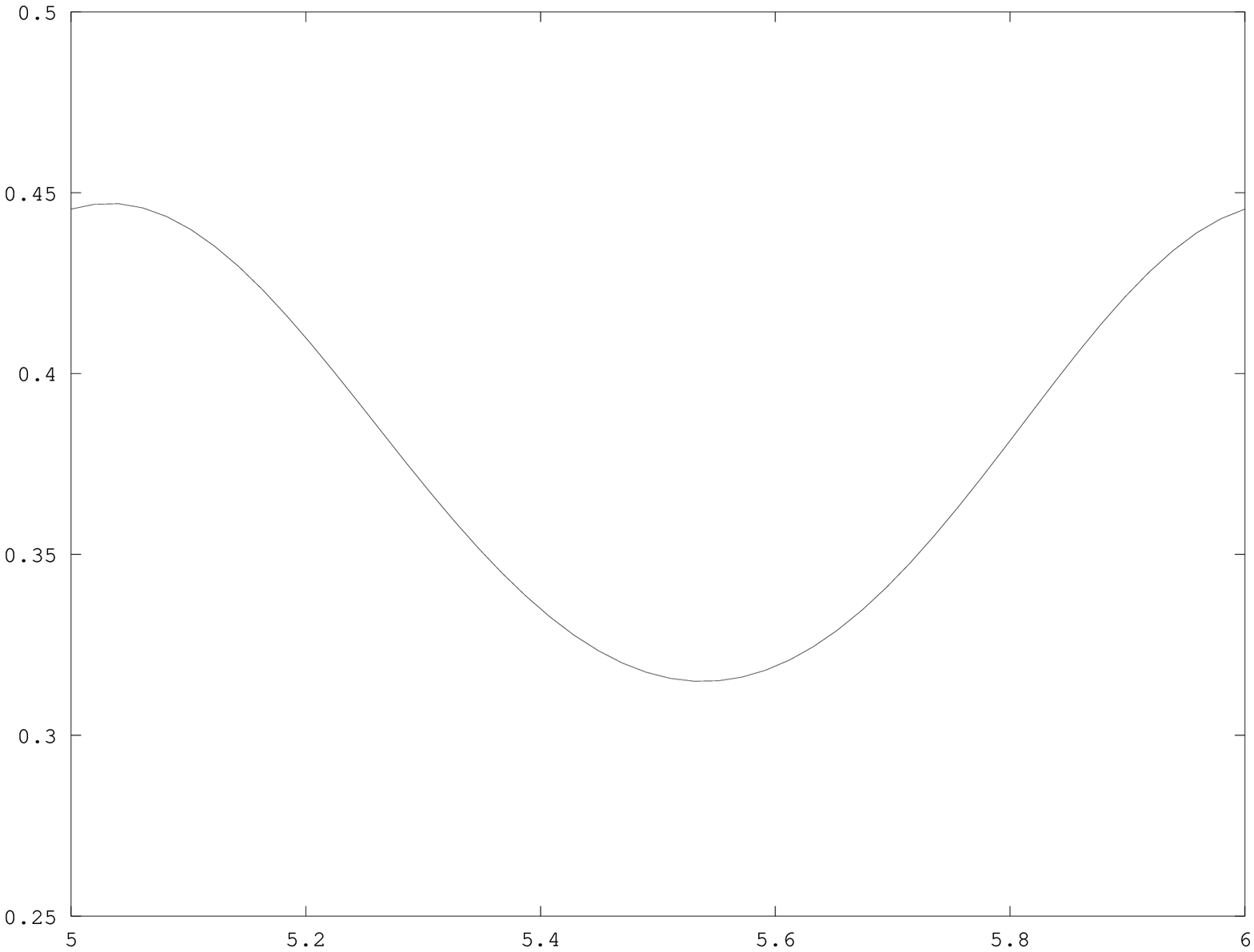}
\end{center}
\caption{Case (ii), the arrival intensity is $\lambda^*(t;20)$. Approximation of the limiting probability
$p_0(t)$ of the empty queue for $t\in[5,6]$.}
\label{fig:4}       % Give a unique label
\end{figure}
\end{center}

\begin{center}
\begin{figure}[!ht]
\begin{center}
  \includegraphics[scale=0.35]{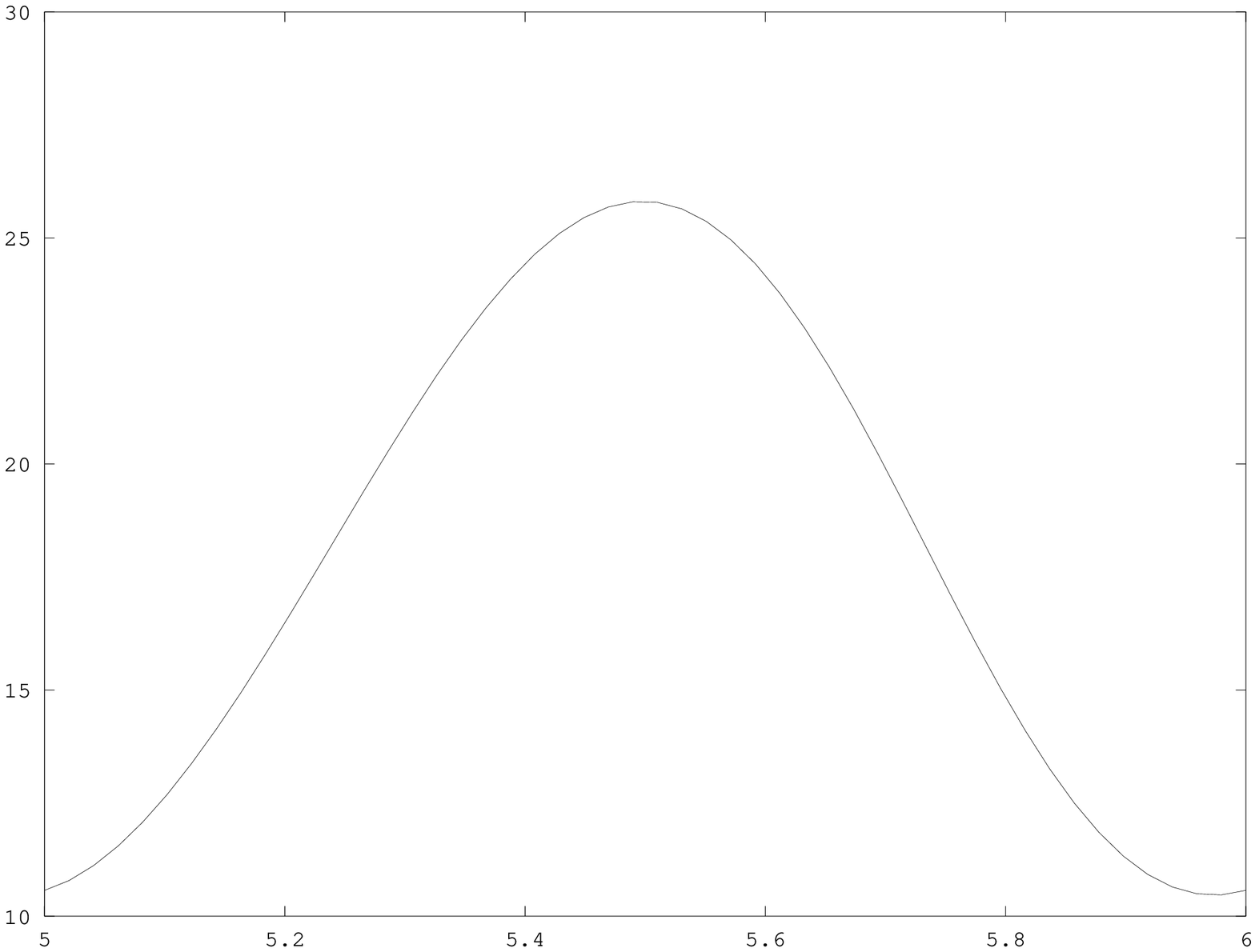}
\end{center}
\caption{Case (ii), the arrival intensity is $\lambda^*(t;50)$. Approximation of the limiting
mean number $\varphi (t)$ of customers in the system for $t\in[5,6]$.}
\label{fig:5}       % Give a unique label
\end{figure}
\end{center}

\begin{center}
\begin{figure}[!ht]
\begin{center}
  \includegraphics[scale=0.35]{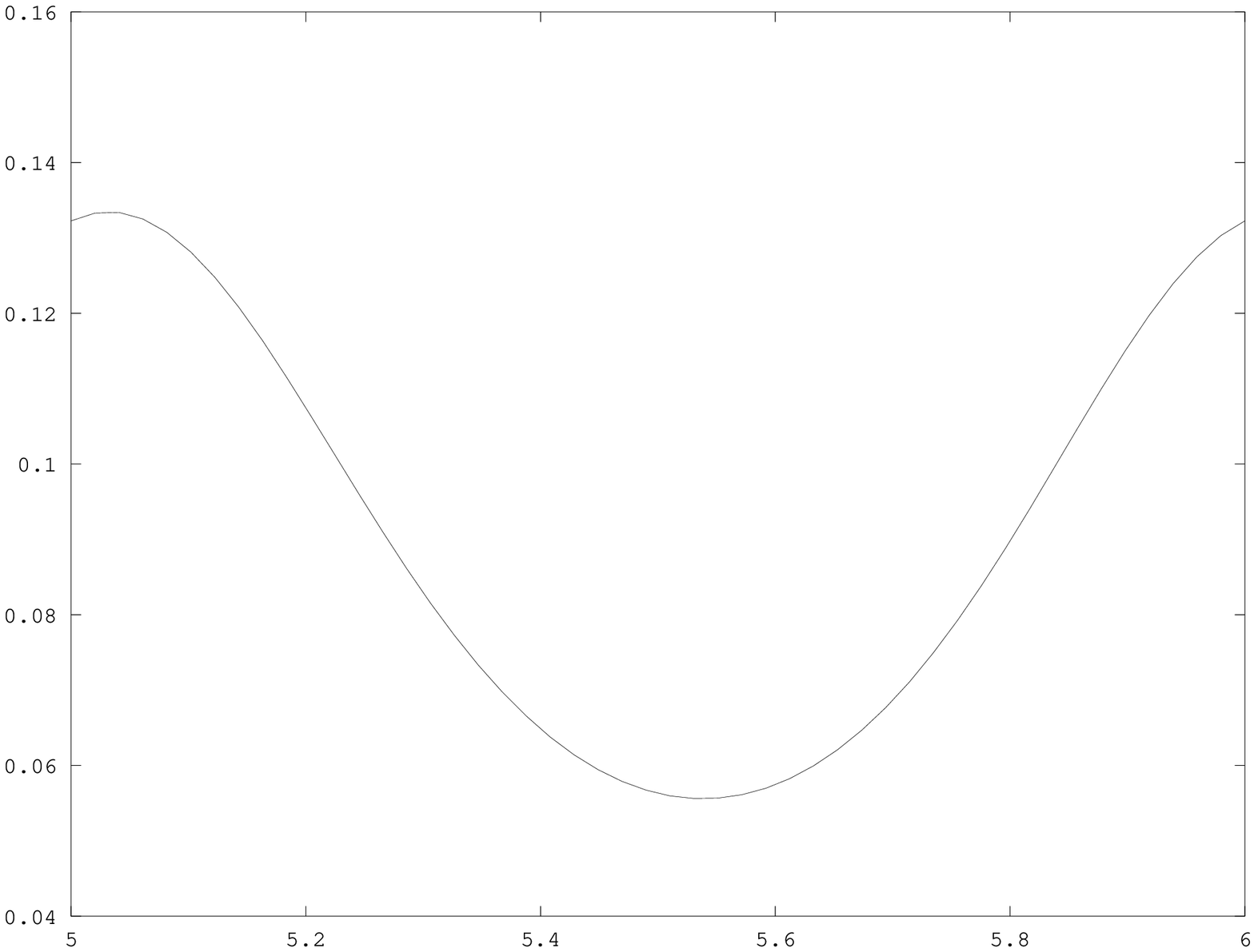}
\end{center}
\caption{Case (ii), the arrival intensity is $\lambda^*(t;50)$. Approximation of the limiting probability
$p_0(t)$ of the empty queue for $t\in[5,6]$.}
\label{fig:6}       % Give a unique label
\end{figure}
\end{center}

%%%%%%%%%% case (iii)
\clearpage

\begin{center}
\begin{figure}[!ht]
\begin{center}
  \includegraphics[scale=0.35]{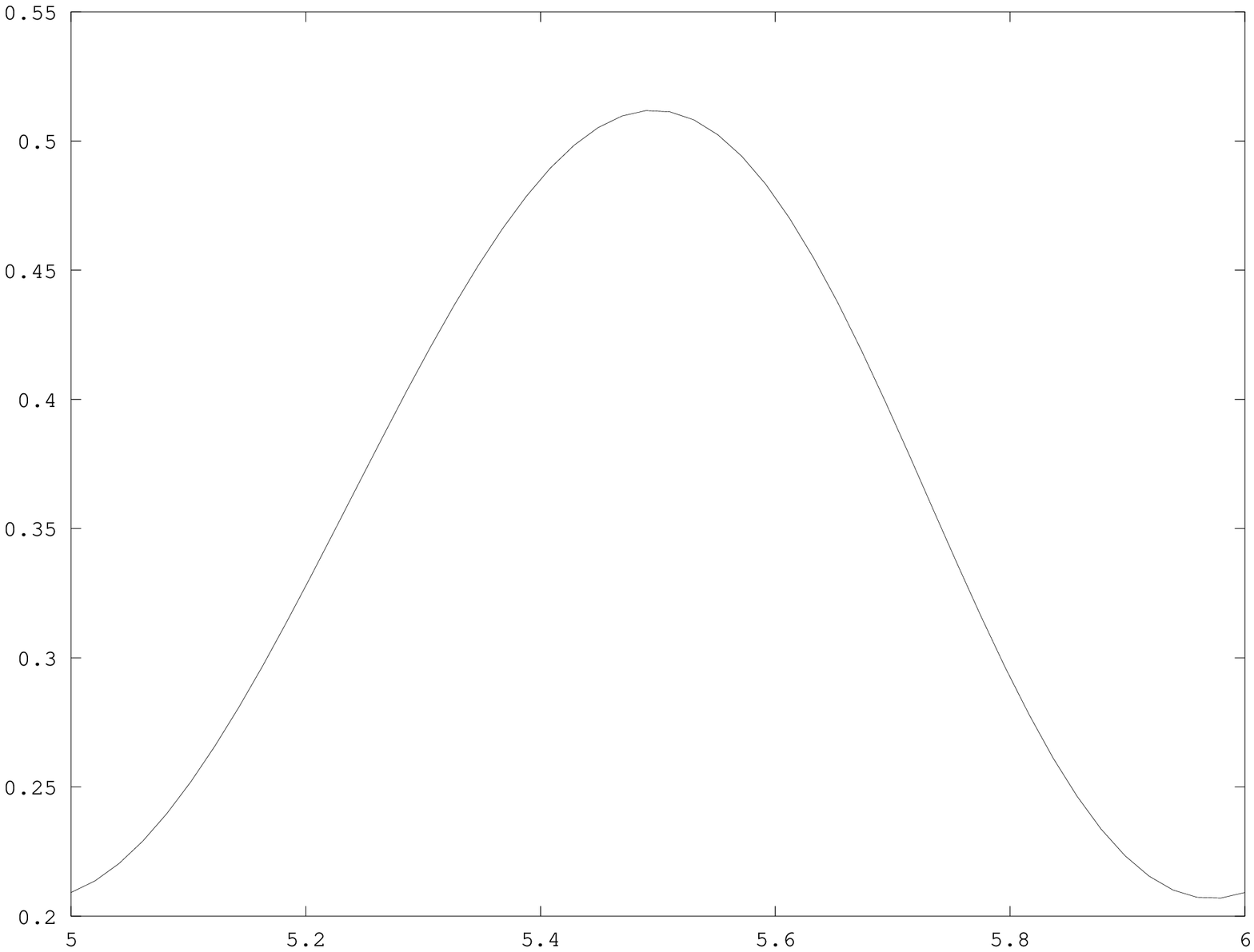}
\end{center}
\caption{Case (iii), the arrival intensity is $\lambda^*(t;10)$. Approximation of the limiting
mean number $\varphi (t)$ of customers in the system for $t\in[5,6]$.}
\label{fig:21}       % Give a unique label
\end{figure}
\end{center}

\begin{center}
\begin{figure}[!ht]
\begin{center}
  \includegraphics[scale=0.35]{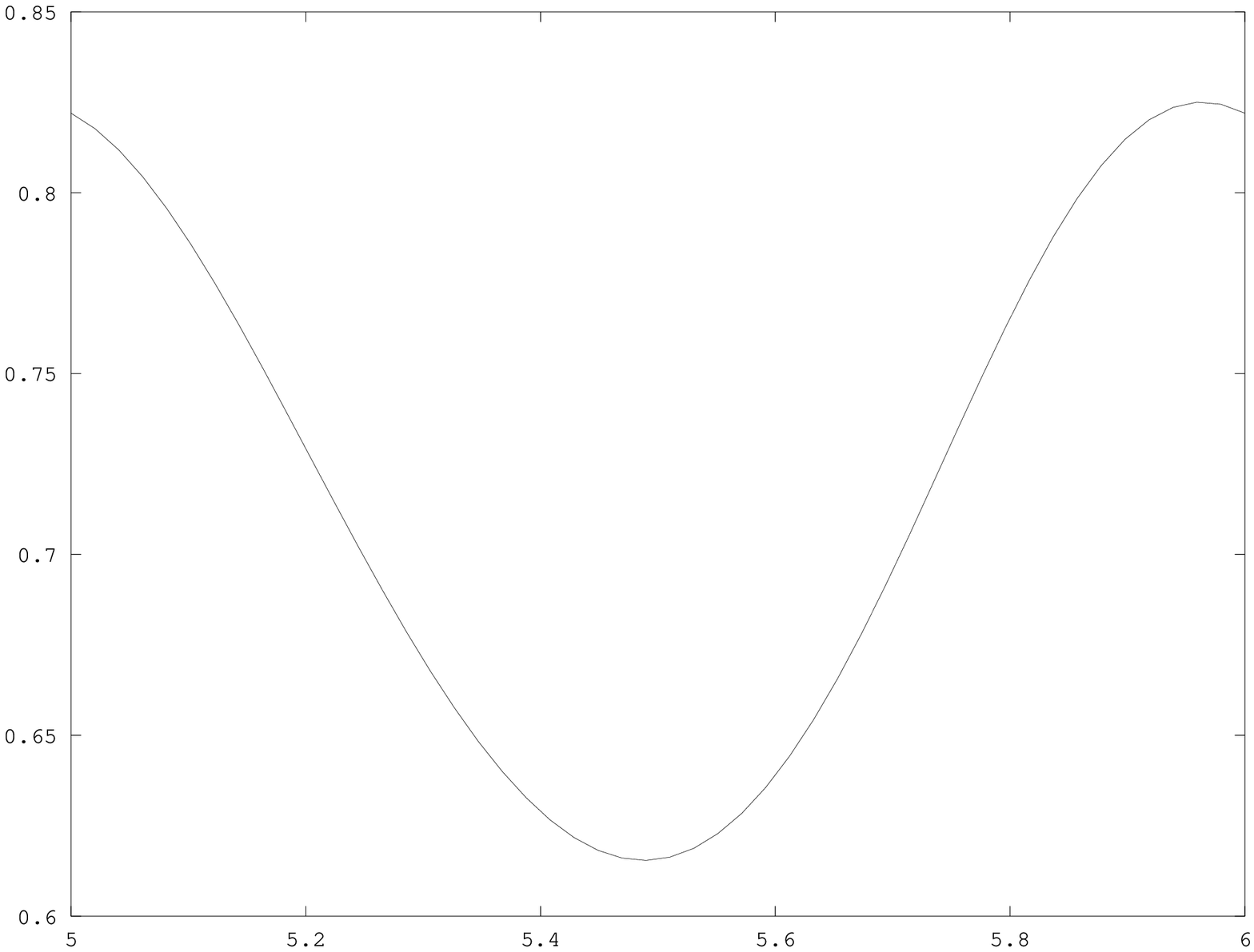}
\end{center}
\caption{Case (iii), the arrival intensity is $\lambda^*(t;10)$. Approximation of the limiting probability
$p_0(t)$ of the empty queue for $t\in[5,6]$.}
\label{fig:22}       % Give a unique label
\end{figure}
\end{center}

\begin{center}
\begin{figure}[!ht]
\begin{center}
  \includegraphics[scale=0.35]{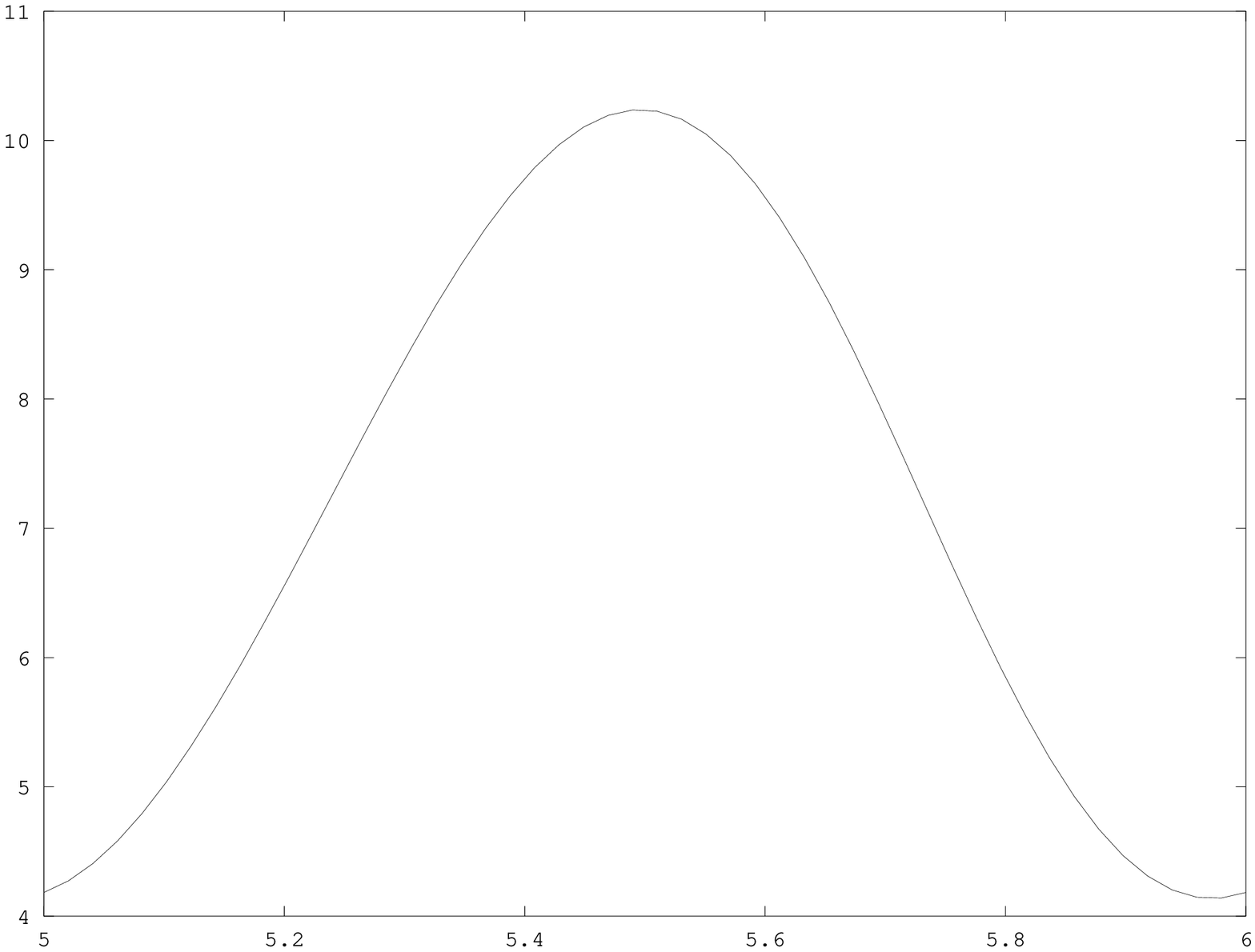}
\end{center}
\caption{Case (iii), the arrival intensity is $\lambda^*(t;20)$. Approximation of the limiting
mean number $\varphi (t)$ of customers in the system for $t\in[5,6]$.}
\label{fig:23}       % Give a unique label
\end{figure}
\end{center}

\begin{center}
\begin{figure}[!ht]
\begin{center}
  \includegraphics[scale=0.35]{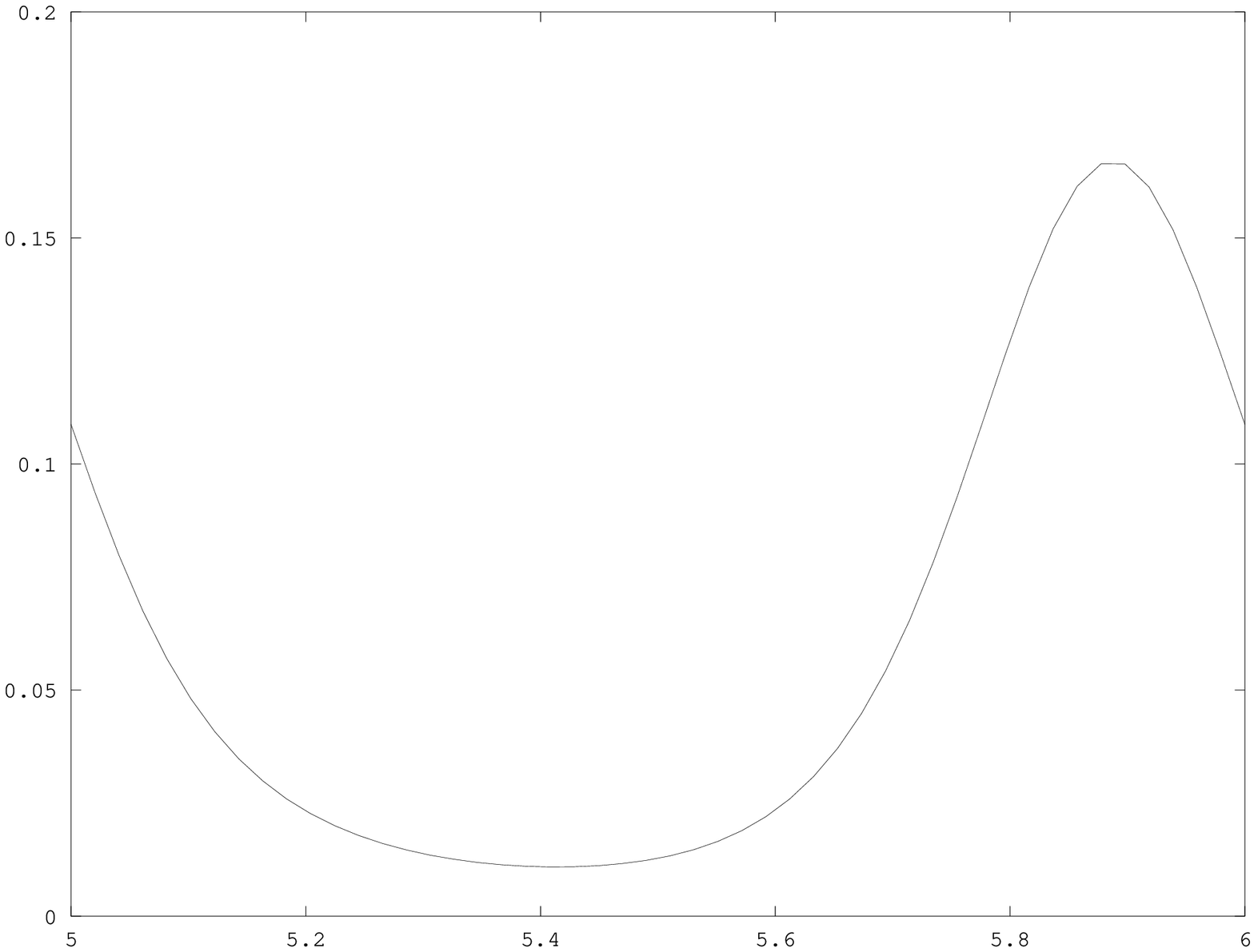}
\end{center}
\caption{Case (iii), the arrival intensity is $\lambda^*(t;20)$. Approximation of the limiting probability
$p_0(t)$ of the empty queue for $t\in[5,6]$.}
\label{fig:24}       % Give a unique label
\end{figure}
\end{center}

\begin{center}
\begin{figure}[!ht]
\begin{center}
  \includegraphics[scale=0.35]{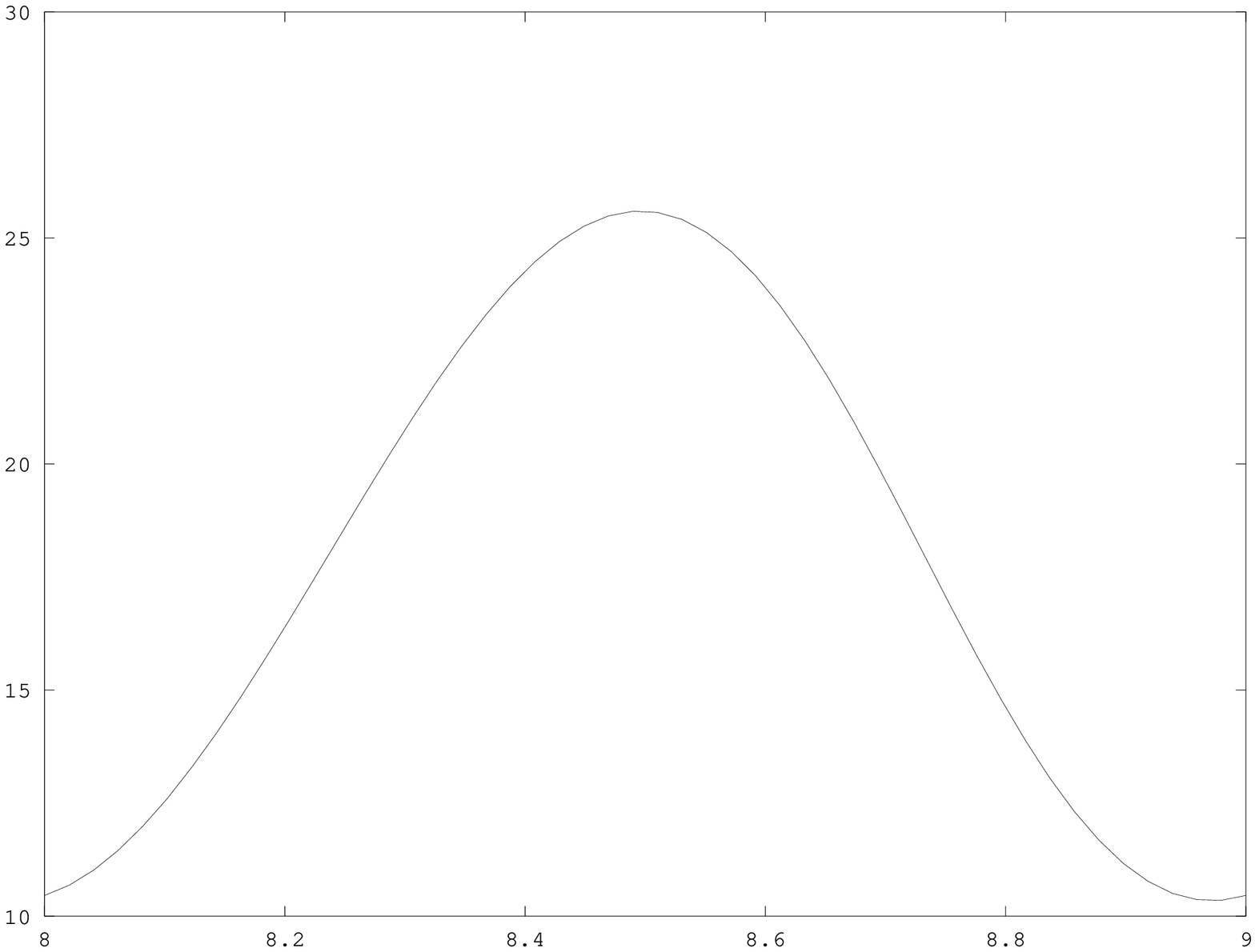}
\end{center}
\caption{Case (iii), the arrival intensity is $\lambda^*(t;50)$. Approximation of the limiting
mean number $\varphi (t)$ of customers in the system for $t\in[8,9]$.}
\label{fig:25}       % Give a unique label
\end{figure}
\end{center}

\begin{center}
\begin{figure}[!ht]
\begin{center}
  \includegraphics[scale=0.35]{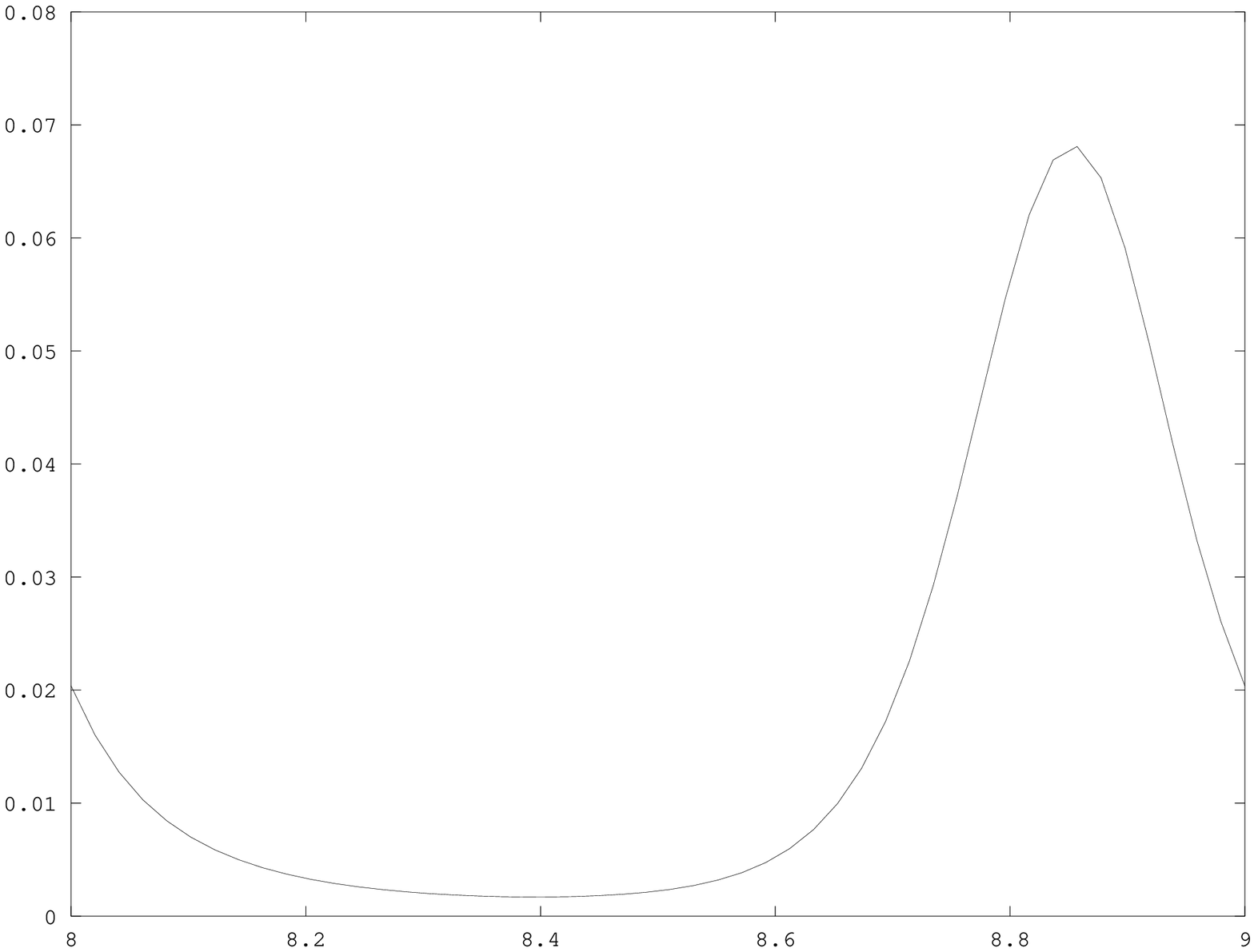}
\end{center}
\caption{Case (iii), the arrival intensity is $\lambda^*(t;50)$. Approximation of the limiting probability
$p_0(t)$ of the empty queue for $t\in[8,9]$.}
\label{fig:26}       % Give a unique label
\end{figure}
\end{center}

%%  case (iv)

\clearpage

\begin{center}
\begin{figure}[!ht]
\begin{center}
  \includegraphics[scale=0.35]{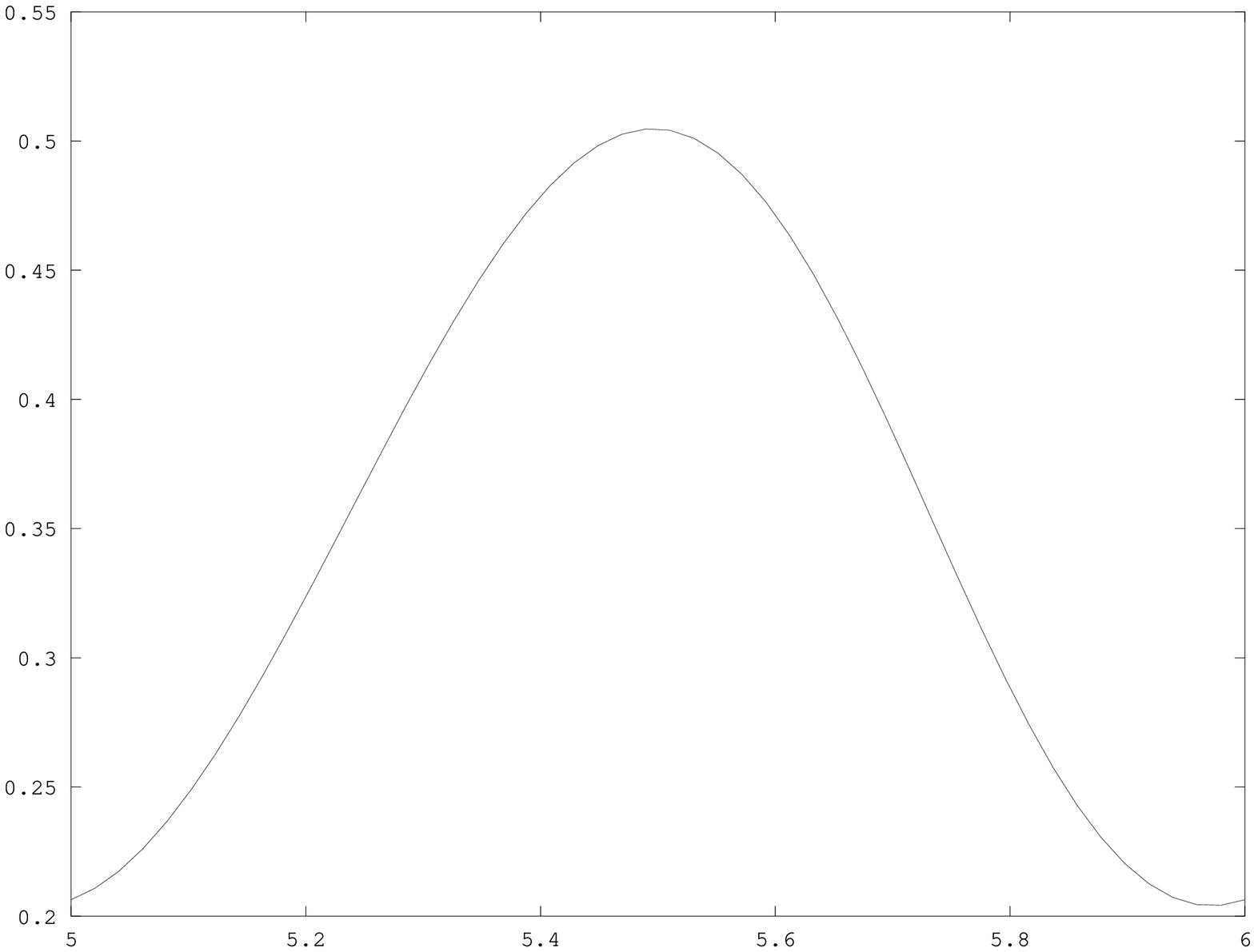}
\end{center}
\caption{Case (iv), the arrival intensity is $\lambda^*(t;10)$. Approximation of the limiting
mean number $\varphi (t)$ of customers in the system for $t\in[5,6]$.}
\label{fig:31}       % Give a unique label
\end{figure}
\end{center}

\begin{center}
\begin{figure}[!ht]
\begin{center}
  \includegraphics[scale=0.35]{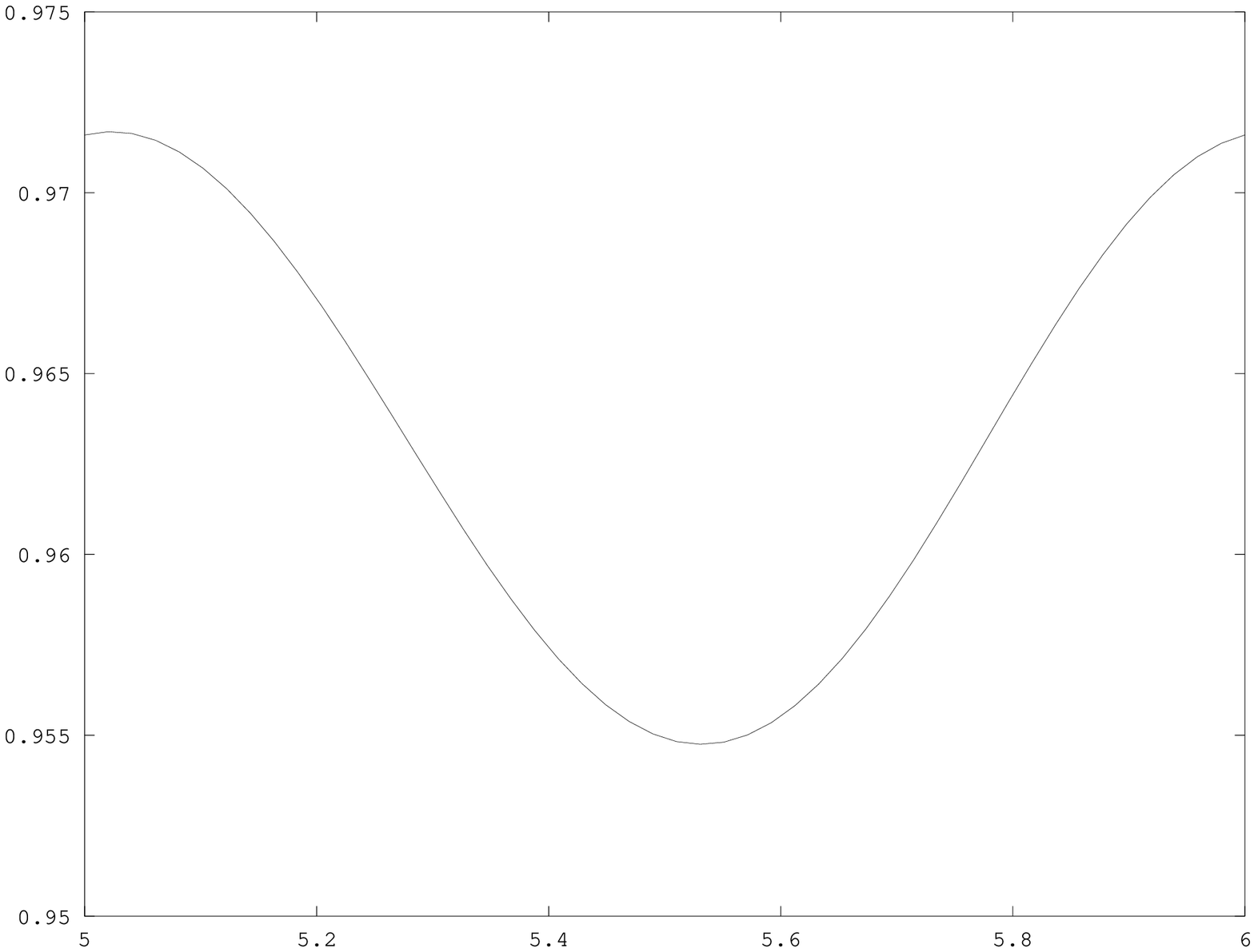}
\end{center}
\caption{Case (iv), the arrival intensity is $\lambda^*(t;10)$. Approximation of the limiting probability
$p_0(t)$ of the empty queue for $t\in[5,6]$.}
\label{fig:32}       % Give a unique label
\end{figure}
\end{center}

\begin{center}
\begin{figure}[!ht]
\begin{center}
  \includegraphics[scale=0.35]{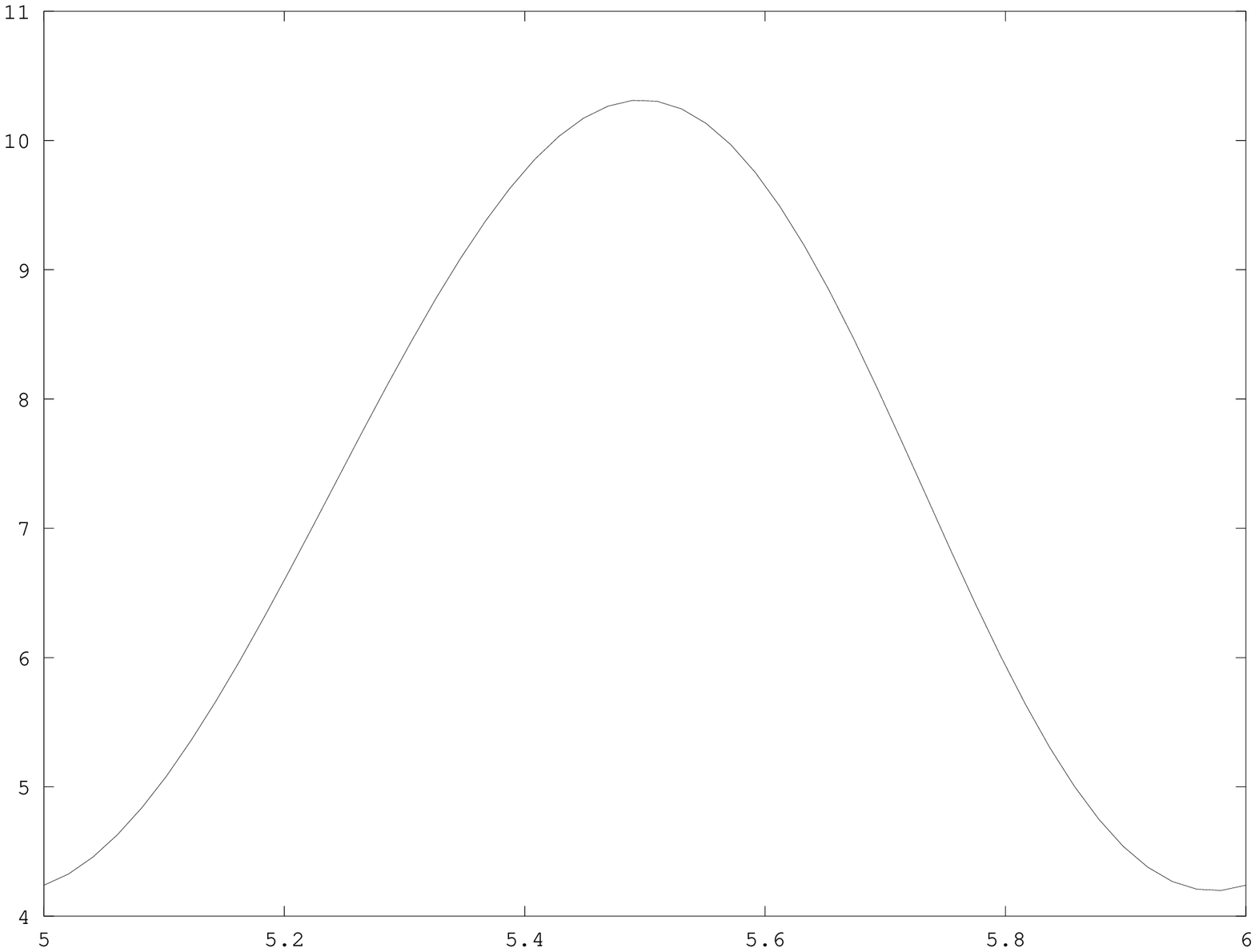}
\end{center}
\caption{Case (iv), the arrival intensity is $\lambda^*(t;20)$. Approximation of the limiting
mean number $\varphi (t)$ of customers in the system for $t\in[5,6]$.}
\label{fig:33}       % Give a unique label
\end{figure}
\end{center}

\begin{center}
\begin{figure}[!ht]
\begin{center}
  \includegraphics[scale=0.35]{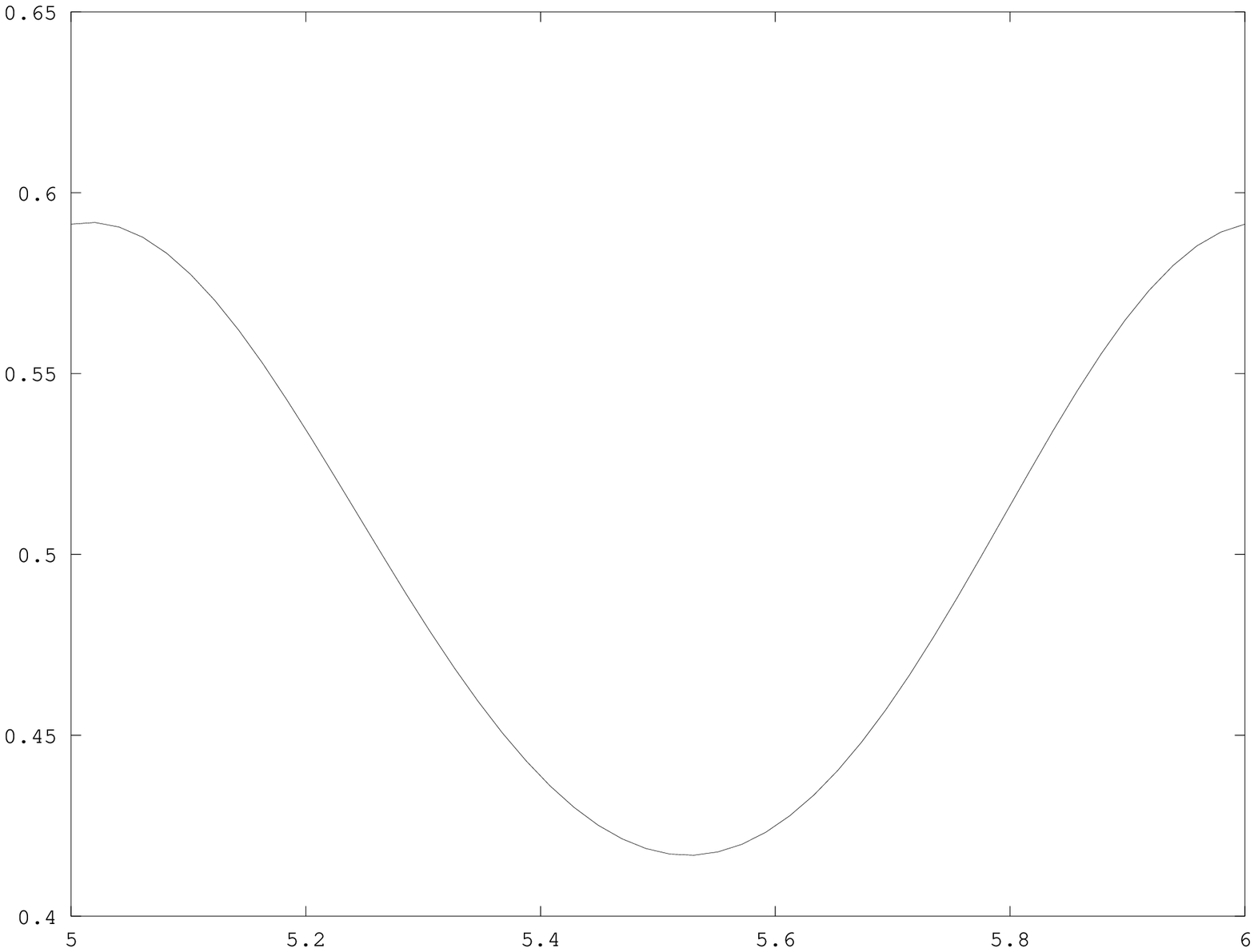}
\end{center}
\caption{Case (iv), the arrival intensity is $\lambda^*(t;20)$. Approximation of the limiting probability
$p_0(t)$ of the empty queue for $t\in[5,6]$.}
\label{fig:34}       % Give a unique label
\end{figure}
\end{center}

\begin{center}
\begin{figure}[!ht]
\begin{center}
  \includegraphics[scale=0.35]{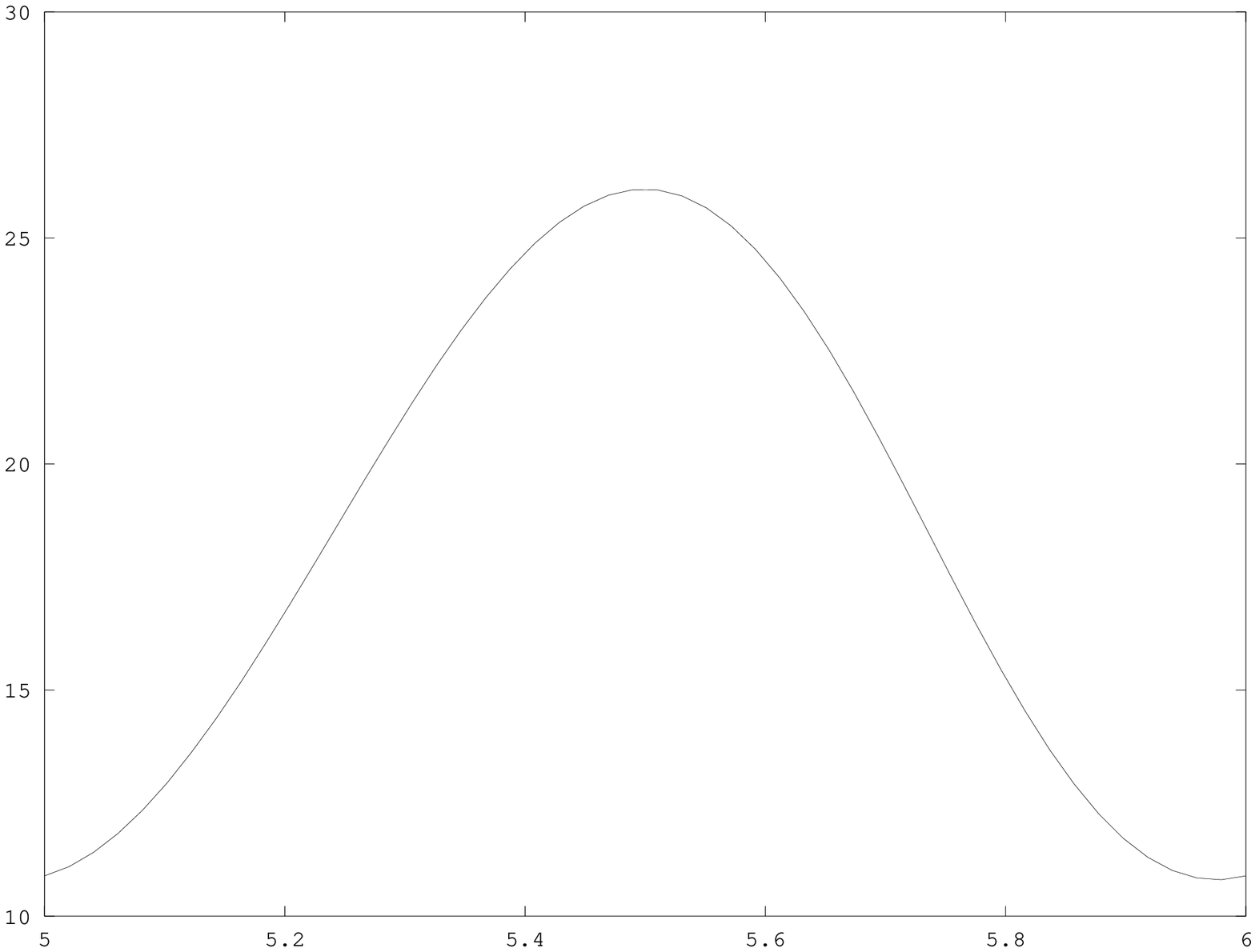}
\end{center}
\caption{Case (iv), the arrival intensity is $\lambda^*(t;50)$. Approximation of the limiting
mean number $\varphi (t)$ of customers in the system for $t\in[5,6]$.}
\label{fig:35}       % Give a unique label
\end{figure}
\end{center}

\begin{center}
\begin{figure}[!ht]
\begin{center}
  \includegraphics[scale=0.35]{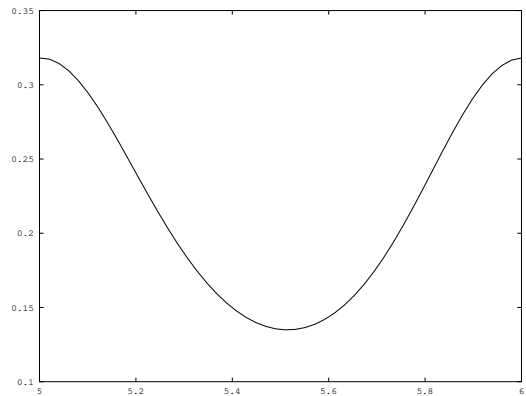}
\end{center}
\caption{Case (iv), the arrival intensity is $\lambda^*(t;50)$. Approximation of the limiting probability
$p_0(t)$ of the empty queue for $t\in[5,6]$.}
\label{fig:36}       % Give a unique label
\end{figure}
\end{center}

\clearpage

\section{Conclusion}

From the presented figures one can see
that the limiting mean number of customers in the system
apparently does not depend on type of the system
i.e. for all four different systems considered
there is numerical evidence that the
limiting means coincide.
With respect to the probability of the empty
queue one observes the clear dependence on
the type of the system.
These numerical evidences show one of the
directions of further research:
explanation these effects from the analytical point of view.
Another direction is the generalization of the
proposed method for other types of inhomogeneous queueing systems.
One of the appealing candidates are queueing systems
 with balking, in which the arrival intensities
 decrease with the growth of the total number of customers
 in the system.
 Another one direction of research
 follows from \cite{Zeifman2016NEW} and is related to the
 optimization of no-wait probabilities in such queueing systems.

%\section*{Acknowledgement}
%This work was supported by Russian Scientific Foundation  (Grant No.
%14-11-00397)??.

\end{document}